\documentclass[a4paper,preprint,fleqn,11pt,authoryear]{elsarticle}
\usepackage{geometry}
\usepackage{fancyhdr}
\usepackage{amsmath}
\usepackage{amssymb}
\usepackage{xcolor}
\usepackage{amssymb}
\usepackage{amsthm}
\usepackage{multirow}
\usepackage{graphicx}
\usepackage[figuresright]{rotating}

\usepackage[linesnumbered,ruled,vlined]{algorithm2e}
\usepackage{booktabs}
\usepackage[font=small,labelfont=bf,labelsep=none]{caption}
\usepackage[figuresright]{rotating}
\usepackage{subfigure}
\usepackage{float}
\usepackage{longtable}
\usepackage{supertabular}

\usepackage{fancybox}
\usepackage{threeparttable}
\usepackage{pdfpages}
\usepackage{makecell}

\usepackage{caption}
\usepackage{makecell}
\usepackage{threeparttable} 
\usepackage{tikz}
\usepackage{rotating}
\usetikzlibrary{shapes.geometric, arrows, positioning}

\usepackage[colorlinks,linkcolor=black,anchorcolor=black,citecolor=black]{hyperref}

\journal{}

\SetKwInput{KwIn}{Input}
\SetKwInput{KwOut}{Output}

\begin{document}
%表格编号Table 1 可以加粗显示
\captionsetup[table]{
	labelsep=newline,%换行
	singlelinecheck=false,%居左
	skip=0.1cm
}

\allowdisplaybreaks[4] %允许公式跨页	
\renewcommand{\thefootnote}{*}
	
\begin{frontmatter}
	
\title{A survey of the orienteering problem: model evolution, algorithmic advances, and future directions}
\author{Songhao Shen$^{a}$, Yufeng Zhou$^{a}$, Qin Lei$^{a}$, Zhibin Wu$^{a,}$\footnote{Corresponding author. \newline \hspace*{1em} E-mail address: zhibinwu@scu.edu.cn (Z.Wu)}}
\address{
	$^{a}$Business School, Sichuan University, Chengdu 610065, China}
%	$^{b}$School of Computing and Information Systems, Singapore Management University, Singapore 178902, Singapore
%}

%\maketitle

\begin{abstract}
The orienteering problem (OP) is a combinatorial optimization problem that seeks a path visiting a subset of locations to maximize collected rewards under a limited resource budget. 
This article presents a systematic PRISMA-based review of OP research published between 2017 and 2025, with a focus on models and methods that have shaped subsequent developments in the field. 
We introduce a component-based taxonomy that decomposes OP variants into time-, path-, node-, structure-, and information-based extensions. 
This framework unifies classical and emerging variants—including stochastic, time-dependent, Dubins, Set, and multi-period OPs—within a single structural perspective. 
We further categorize solution approaches into exact algorithms, heuristics and metaheuristics, and learning-based methods, with particular emphasis on matheuristics and recent advances in artificial intelligence, especially reinforcement learning and neural networks, which enhance scalability in large-scale and information-rich settings. 
Building on this unified view, we discuss how different components affect computational complexity and polyhedral properties and identify open challenges related to robustness, sustainability, and AI integration. 
The survey thus provides both a consolidated reference for existing OP research and a structured agenda for future theoretical and applied work.
\end{abstract}


\end{frontmatter}

\section{Introduction}
\label{sec:1}

\begin{figure*}[htbp]
	%[!t]
	\centering
	\subfigure[Traveling Salesman Problem]{
		%		\label{fig:operatora} %% label for first subfigure
		\includegraphics[width=0.25\linewidth]{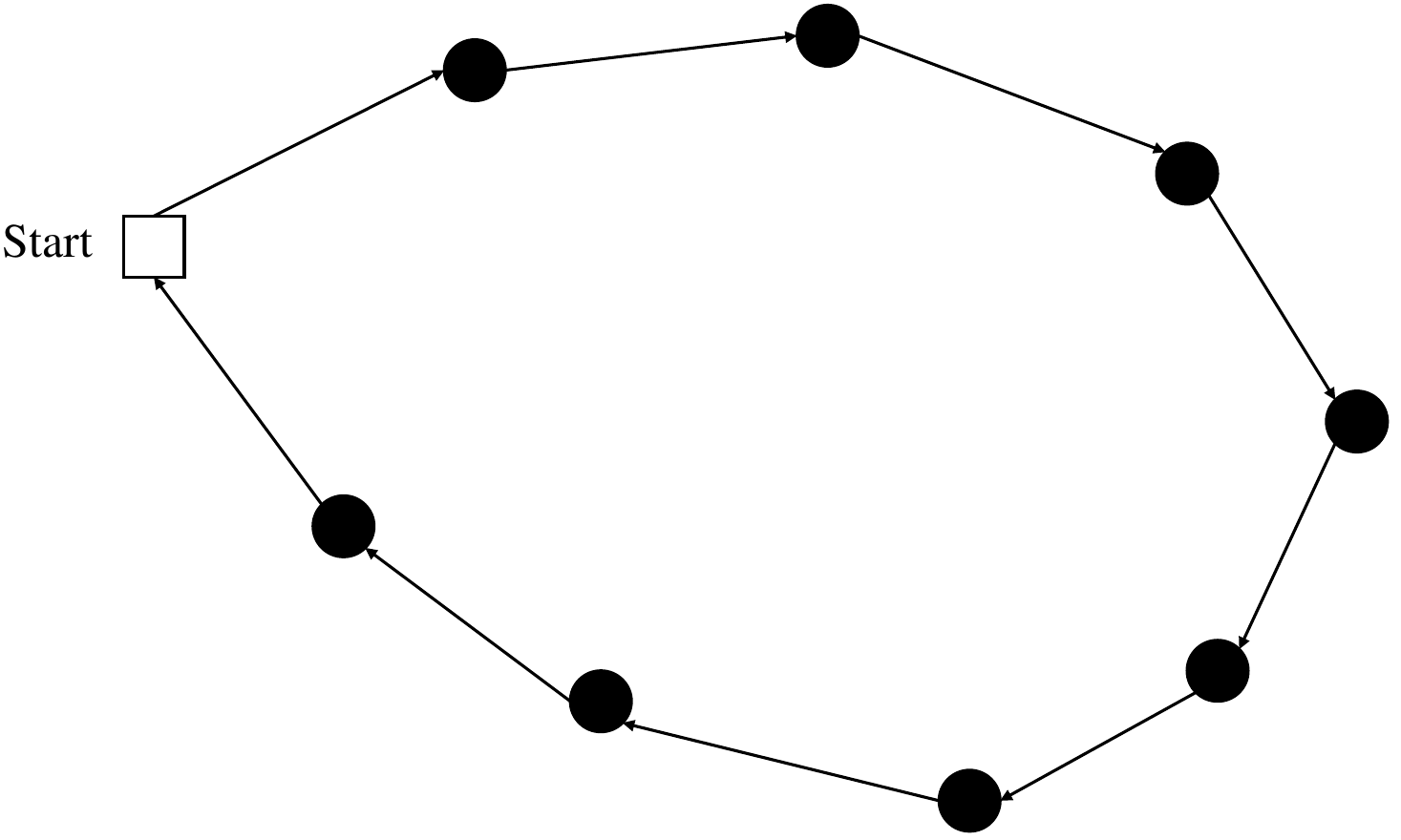}}
	\subfigure[Vehicle Routing Problem]{
		%		\label{fig:operatorc} %% label for first subfigure
		\includegraphics[width=0.25\linewidth]{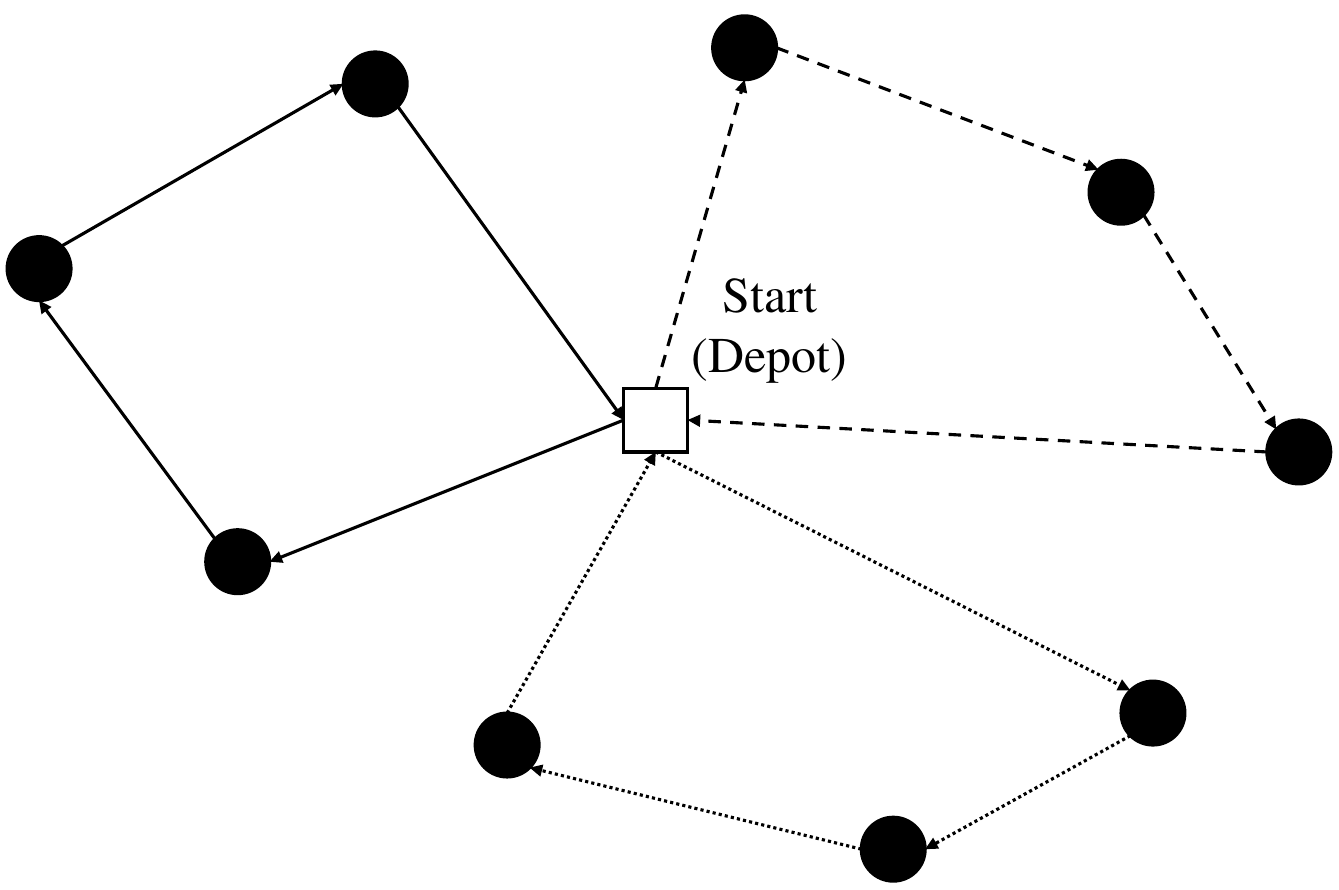}}\\
	\subfigure[Orienteering Problem]{
		%		\label{fig:operatora} %% label for first subfigure
		\includegraphics[width=0.25\linewidth]{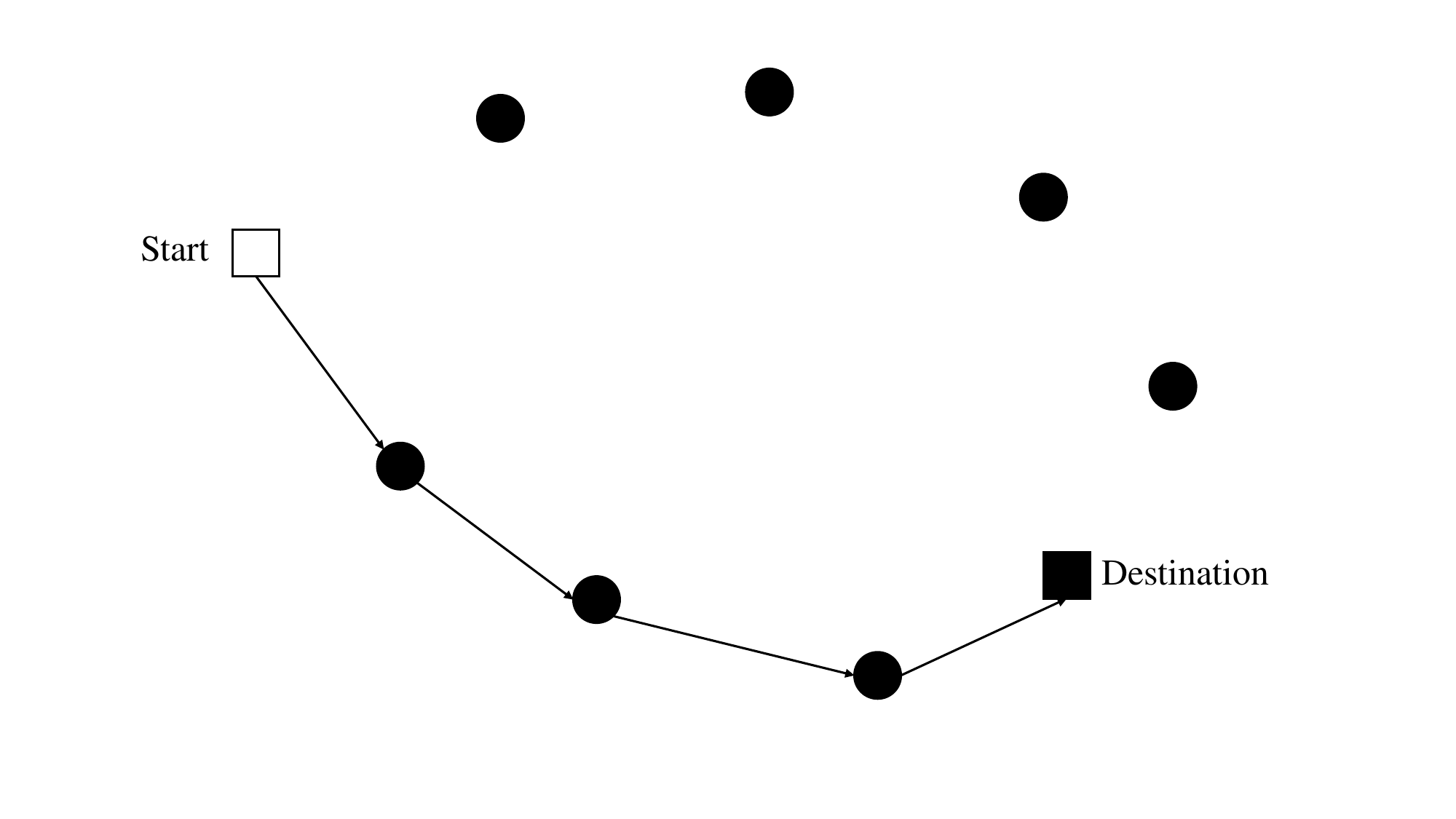}}
	\subfigure[Team Orienteering Problem]{
		%		\label{fig:operatorc} %% label for first subfigure
		\includegraphics[width=0.25\linewidth]{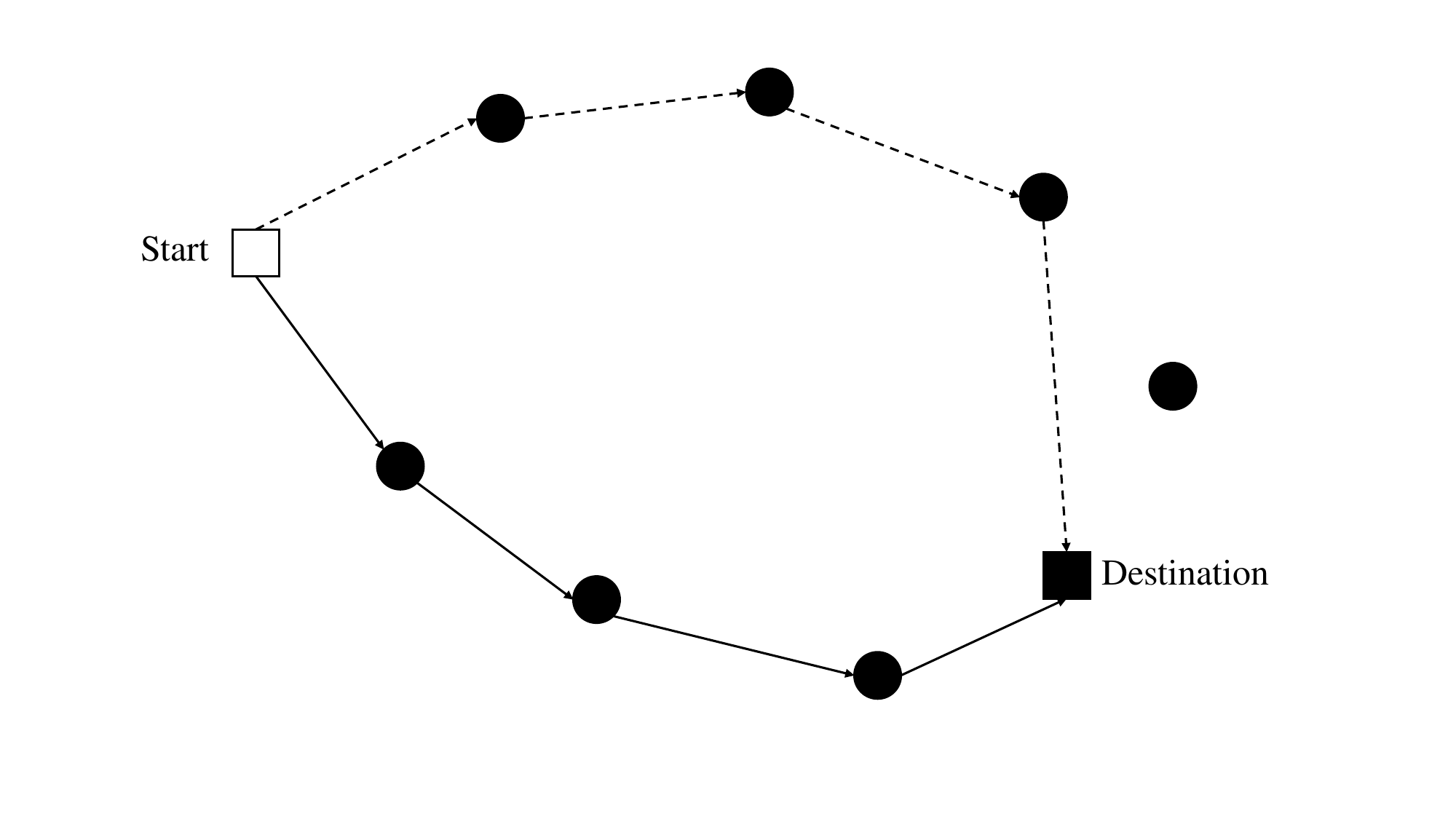}}
	\caption{\quad Example of OP and its related combinatorial optimization problems.}
	\label{fig:fig1} %% label for entire figure
\end{figure*}

The term ``orienteering problem'' (OP) was first introduced by \cite{Golden1987}. Since then, the OP and its variants have advanced substantially in both theory and practice. The OP is a classical combinatorial optimization problem derived from the well-known traveling salesman problem (TSP). Unlike the TSP or vehicle routing problem (VRP), which aim to minimize the cost of serving all locations, the OP seeks a path visiting only a subset of locations so as to maximize total collected reward under a resource budget (e.g., time or distance). This selectivity makes the OP a natural framework for resource-constrained applications such as urban logistics, tourism, and unmanned aerial vehicle (UAV) path planning (see Fig.~\ref{fig:fig1}).

\begin{figure*}[!tbp]
	\centering
	\includegraphics[width=0.75\linewidth]{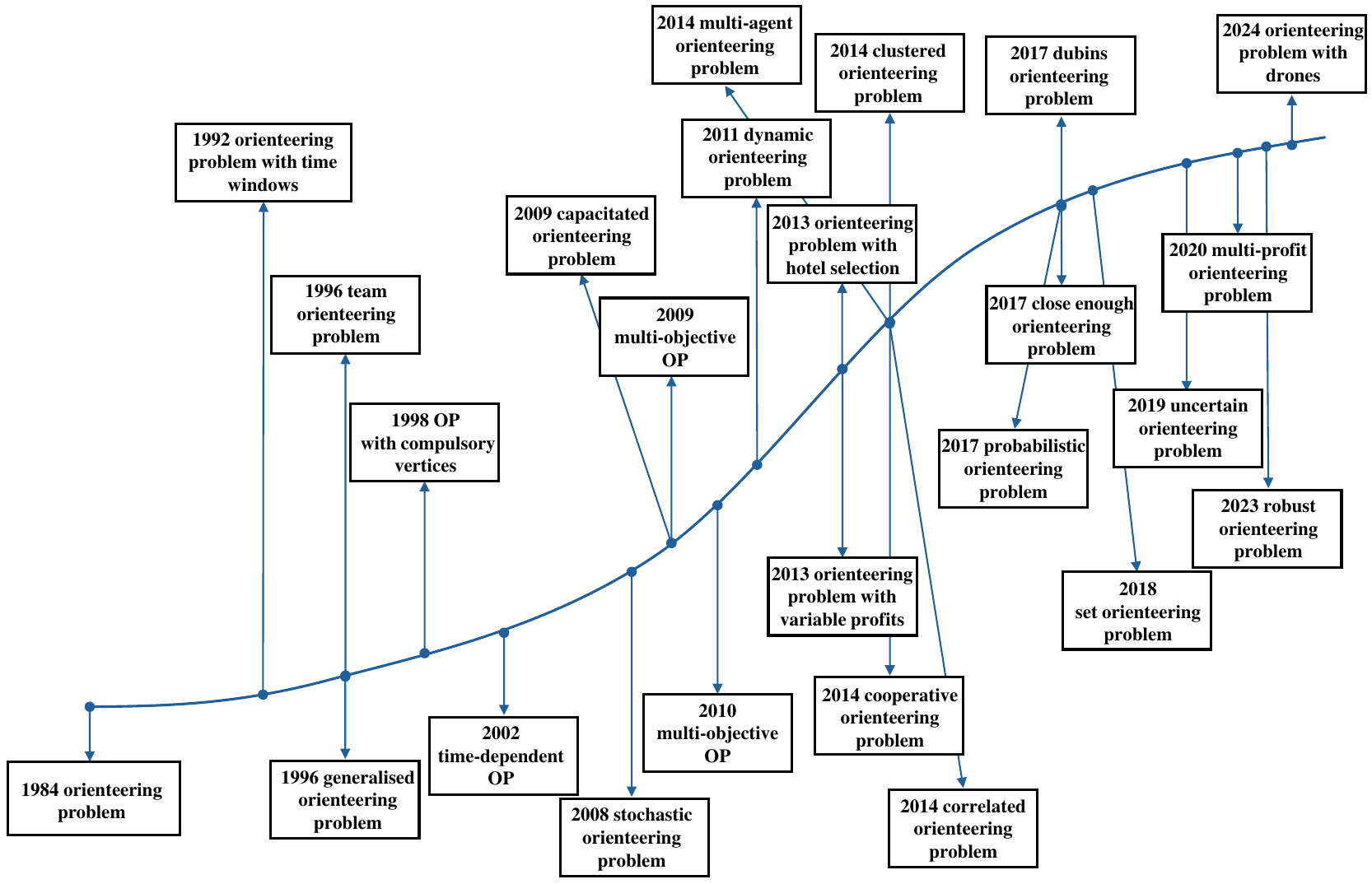}
	\caption{\quad Timeline of concept introductions in extensions of the Orienteering Problems.}
	\label{fig_timeline}
\end{figure*}

Over recent decades, OP research has evolved from deterministic, single-vehicle models to a broad family of variants that incorporate rich temporal, spatial, and informational constraints (Fig.~\ref{fig_timeline}). Foundational components—such as time windows, team and capacity constraints—have been progressively combined with emerging stochastic, set-based, and dynamic features to capture complex real-world environments. These developments have been accompanied by a wide spectrum of solution techniques, ranging from exact algorithms and matheuristics to reinforcement learning and other AI-driven approaches.

Earlier surveys have provided comprehensive overviews of OP variants and algorithms up to around 2017 (e.g., \cite{vansteenwegen2011,gunawan2016}). Since then, however, the literature has expanded rapidly in several directions: (i) the introduction of new structural variants such as Set, Dubins, and drone OPs; (ii) the systematic integration of uncertainty and multi-period decision making; and (iii) the adoption of learning-based methods that exploit large-scale and information-rich settings. Despite this progress, there is still a lack of a unified structural framework that connects these recent models and methods and clarifies how their underlying components affect computational properties and algorithm design.

This article fills that gap by conducting a systematic review of OP research from 2017 to 2025 and by organizing recent developments within a component-based structural taxonomy. Using a PRISMA-style protocol, we identify and analyze 112 peer-reviewed studies that represent the state of the art in OP modeling and solution methodologies. Building on this corpus, we address the following research questions:
\begin{itemize}
	\item RQ1: How have the theoretical models of the OP evolved from 2017 to 2025, and what are the emerging foundational components?
	\item RQ2: What are the dominant solution strategies (exact, heuristic, and learning-based) applied to these variants?
	\item RQ3: What are the current research foci, gaps, and future trends in the OP domain?
\end{itemize}

The main contributions to the orienteering literature are as follows.

\begin{enumerate}
	\item We propose a component-based taxonomy—comprising time-, path-, node-, structure-, and information-based extensions—that unifies classical and emerging OP variants (such as set, Dubins, and drone OPs) within a single framework.
	
	\item We analyze how each component alters the OP’s computational complexity and polyhedral properties. This structural perspective connects OP variants to canonical problem families and explains the affinity of specific variants for particular exact formulations or approximation paradigms.

	\item We synthesize recent advances in exact, heuristic, and learning-based methods, emphasizing matheuristics and reinforcement learning, and outline a research agenda focused on AI integration, scalability, robustness under uncertainty, and sustainability-oriented applications.
\end{enumerate}

\section{Research Design and Component-based Taxonomy}
\label{sec:methodology}

To ensure transparency and rigor, this review follows the PRISMA (Preferred Reporting Items for Systematic Reviews and Meta-Analyses) guidelines. 
Building on the research questions articulated in Section~\ref{sec:1}, this section explains how we designed the review to answer them by detailing the search strategy, selection criteria, and screening process.

\subsection{Search Strategy and Selection Criteria}

The literature search was conducted in the \textit{Web of Science Core Collection} (SCI-Expanded and SSCI), which indexes high-impact journals in operations research and related fields.
We focused on articles published between 2017 and 2025 whose topic is related to the ``orienteering problem'' and its recognized variants.

Studies were included if they satisfied all of the following criteria:
\begin{itemize}
	\item \textbf{Topic relevance:} the title, abstract, or keywords explicitly refer to the orienteering problem or a well-established OP variant.
	\item \textbf{Publication type and quality:} peer-reviewed journal article written in English.
	\item \textbf{Methodological contribution:} the paper proposes a new model, algorithm, or a substantial methodological improvement beyond straightforward parameter tuning or minor implementation changes.
\end{itemize}

Studies were excluded if any of the following conditions held:
\begin{itemize}
	\item \textbf{Application-only focus:} reports that use OP-based models purely as a tool for a specific case study without methodological novelty.
	\item \textbf{Isolated or overly specific variants:} single-paper special cases with niche constraints that have not been adopted or generalized in subsequent literature.
	\item \textbf{Outside optimization scope:} works without an explicit optimization model or algorithmic component, even if they are loosely related to routing or planning.
\end{itemize}

To balance established works and cutting-edge developments, we adopted a dynamic citation threshold strategy.
Older papers were required to meet a higher minimum citation count, whereas more recent publications could be retained with fewer citations, especially if they appeared in leading (Q1--Q2) journals.
This approach acknowledges that citation accumulation takes time and prevents very recent yet influential contributions from being discarded solely due to low citation counts.

\begin{figure}[!tbp]
	\centering
	\includegraphics[width=0.7\linewidth]{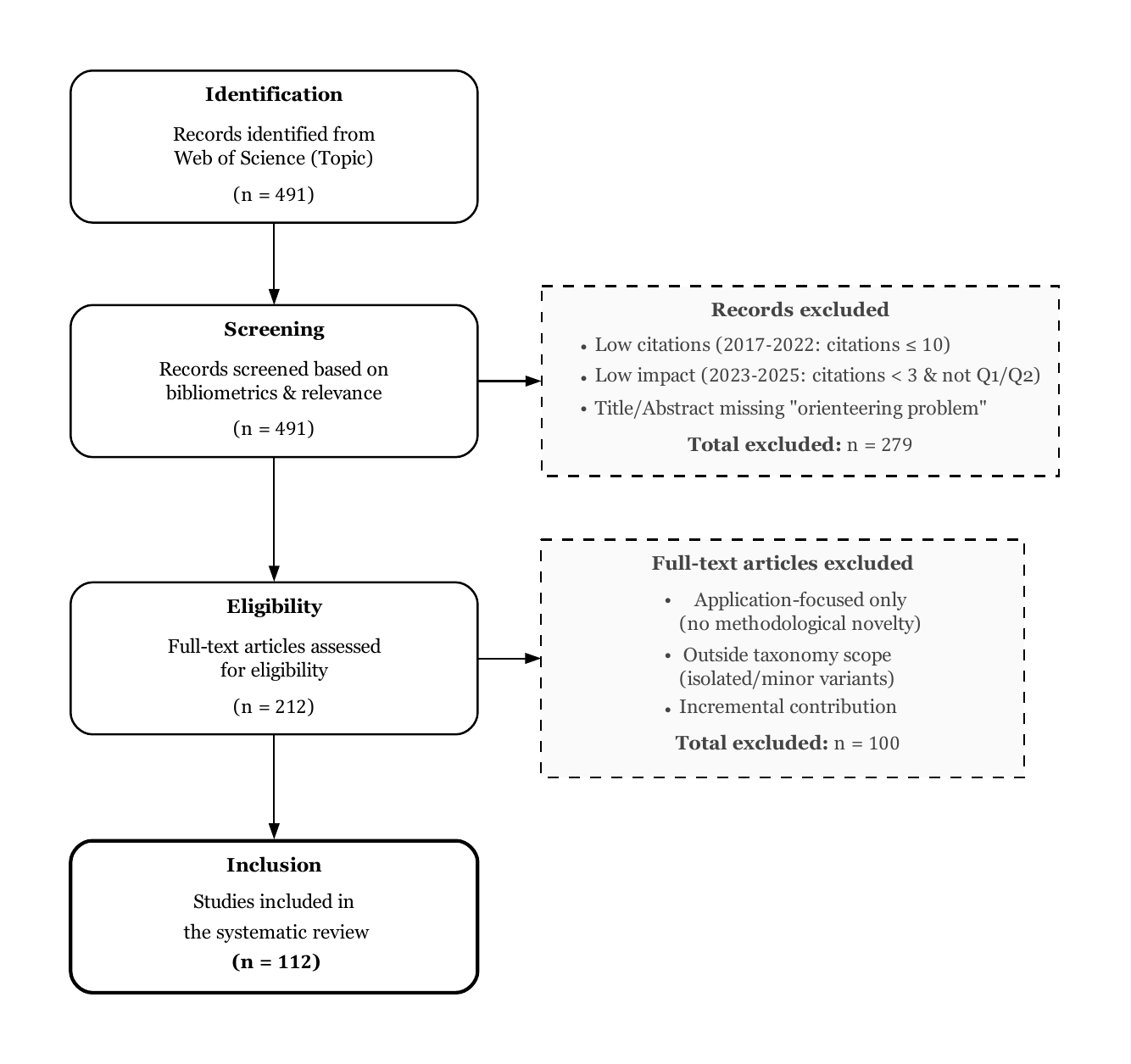} 
	\caption{\quad PRISMA flow diagram of the literature selection process.}
	\label{fig:prisma}
\end{figure}

\subsection{Screening Process}

Study selection followed the multi-stage screening process summarized in Fig.~\ref{fig:prisma}.
The initial keyword search identified 491 records in the Web of Science Core Collection.
After removing duplicates, we screened titles and abstracts against the topic relevance and publication-type criteria described in the previous subsection, discarding clearly irrelevant or non-journal items.

In the next stage, full texts were assessed using the methodological and bibliometric criteria, including the exclusion of application-only articles and isolated model variants.
This ensured that the final sample concentrates on models and algorithms that have had, or are likely to have, a lasting impact on the OP literature.
In addition, citation snowballing was used to identify influential works that were not captured by the database search but clearly satisfied all inclusion criteria.
While not exhaustive, the resulting set of studies provides a comprehensive and representative coverage of recent research on orienteering problems.

\subsection{Dataset Description}
OP research is concentrated in high-impact journals, with the \textit{European Journal of Operational Research}, \textit{Computers \& Operations Research}, and \textit{Applied Soft Computing} being the primary outlets.

\begin{figure*}[!tbp]
	\centering
	\includegraphics[width=0.85\linewidth]{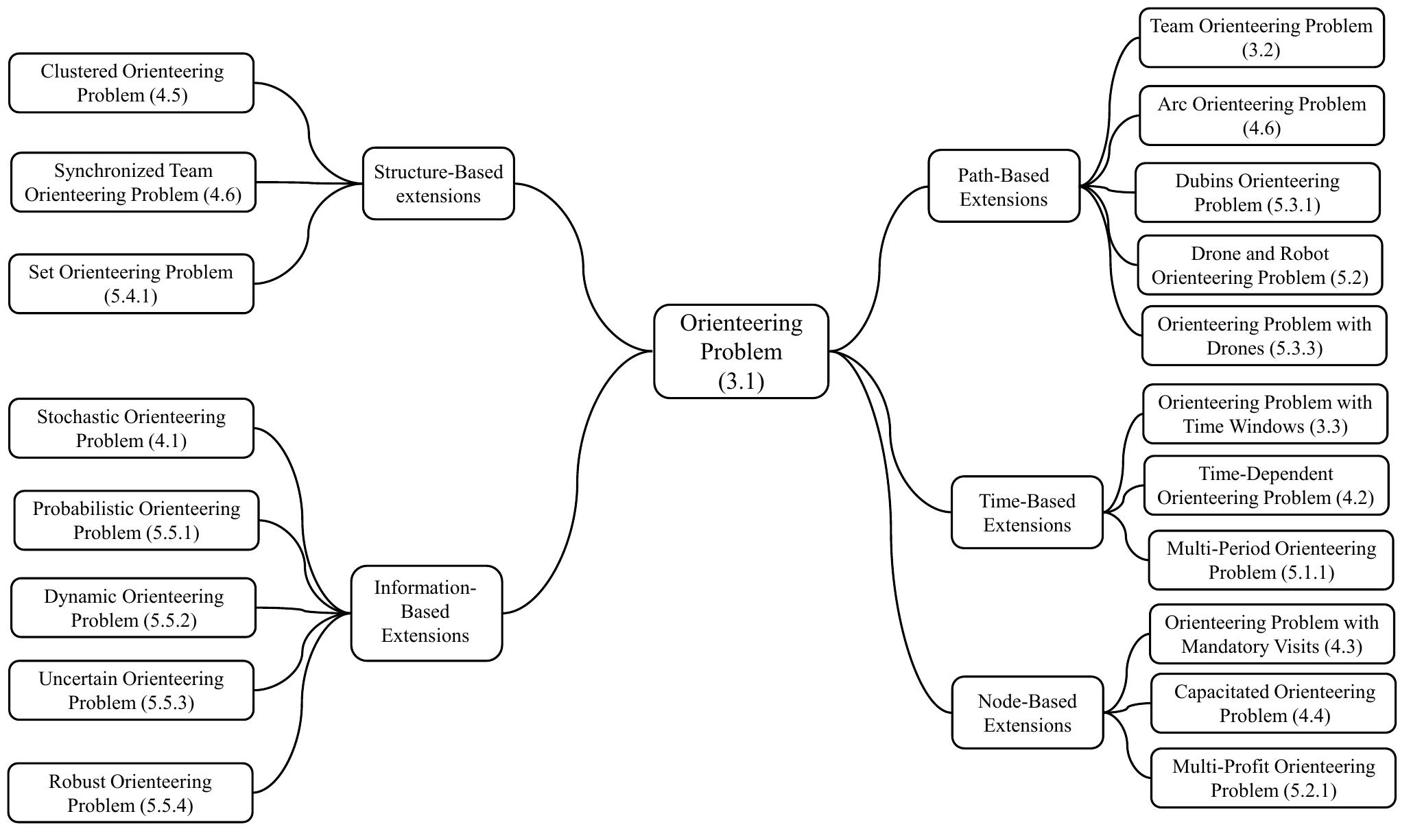}
	\caption{\quad Extensions of the Orienteering Problem.}
	\label{fig_fig3}
\end{figure*}

Fig. \ref{fig_fig3} outlines the primary OP extensions, with bracketed numbers indicating the relevant sections.
These variants function independently but often intersect.
This allows for complex, context-sensitive models tailored to applications like logistics, disaster response, and urban planning.

%\begin{table}[htbp]
%	\centering
%	\caption{Overview of journals publishing research on the Orienteering Problem..}
%	\scriptsize % 1. 使用小号字体 (如果还觉得大，可以改成 \tiny)
%	\setlength{\tabcolsep}{1.5pt} % 2. 极大地减小列间距 (默认是 6pt)
%	\renewcommand{\arraystretch}{0.75}
%	\begin{tabular}{lc}
%		\toprule
%		\multicolumn{1}{c}{Journal } & \multicolumn{1}{c}{Number of Articles} \\
%		\midrule
%		European Journal of Operational Research & 73 \\
%		Computers \& Operations Research & 64 \\
%		Applied Soft Computing & 25 \\
%		Computers \& Industrial Engineering & 24 \\			
%		Expert Systems With Applications & 21 \\
%		Transportation Science & 18 \\
%		Networks & 17 \\
%		Mathematics & 16 \\
%		Annals of Operations Research & 14 \\
%		IEEE Robotics and Automation Letters & 13 \\
%		Journal of The Operational Research Society & 12 \\					
%		Transportation Research Part E: Logistics and Transportation Review & 12 \\
%		Journal of Heuristics & 10 \\
%		IEEE Access & 9 \\
%		Applied Sciences & 9 \\
%		International Transactions In Operational Research & 8 \\
%		Omega: International Journal of Management Science & 8 \\
%		Transportation Research Part B: Methodological & 8 \\
%		INFORMS Journal on Computing & 7 \\
%		International Journal of Production Research & 7 \\
%		IEEE Transactions on Automation Science and Engineering & 6 \\
%		IEEE Transactions on Intelligent Transportation Systems & 6 \\
%		IEEE Transactions on Robotics & 6 \\
%		Military Operations Research & 6 \\
%		OR Spectrum & 6 \\
%		Tourism Management & 6 \\			
%		\bottomrule
%	\end{tabular}%
%	\label{tab:journals}%
%\end{table}%

We classify these extensions into five categories:

(1) Time-based extensions: Including time windows, time-dependent, and multi-period OPs, these primarily involve temporal factors.

(2) Path-based extensions: Covering arc, multi-path, and drone-related OPs, these focus on physical path characteristics.

(3) Node-based extensions: Involving mandatory visits, capacity constraints, and Steiner OPs, these address specific node requirements.

(4) Structure-based extensions: Comprising set, clustered, and synchronized team OPs, these address structural correlations, such as group dependencies or simultaneous visits.

(5) Information-acquisition-based extensions: Including stochastic, dynamic, and robust OPs, these focus on information quality, such as dynamic updates or probabilistic data.

\begin{figure*}[!thbp]
	%[!t]
	\centering
	\subfigure[Article distribution by theme and year.]{
		\label{fig:operatora}
		\includegraphics[width=0.45\linewidth]{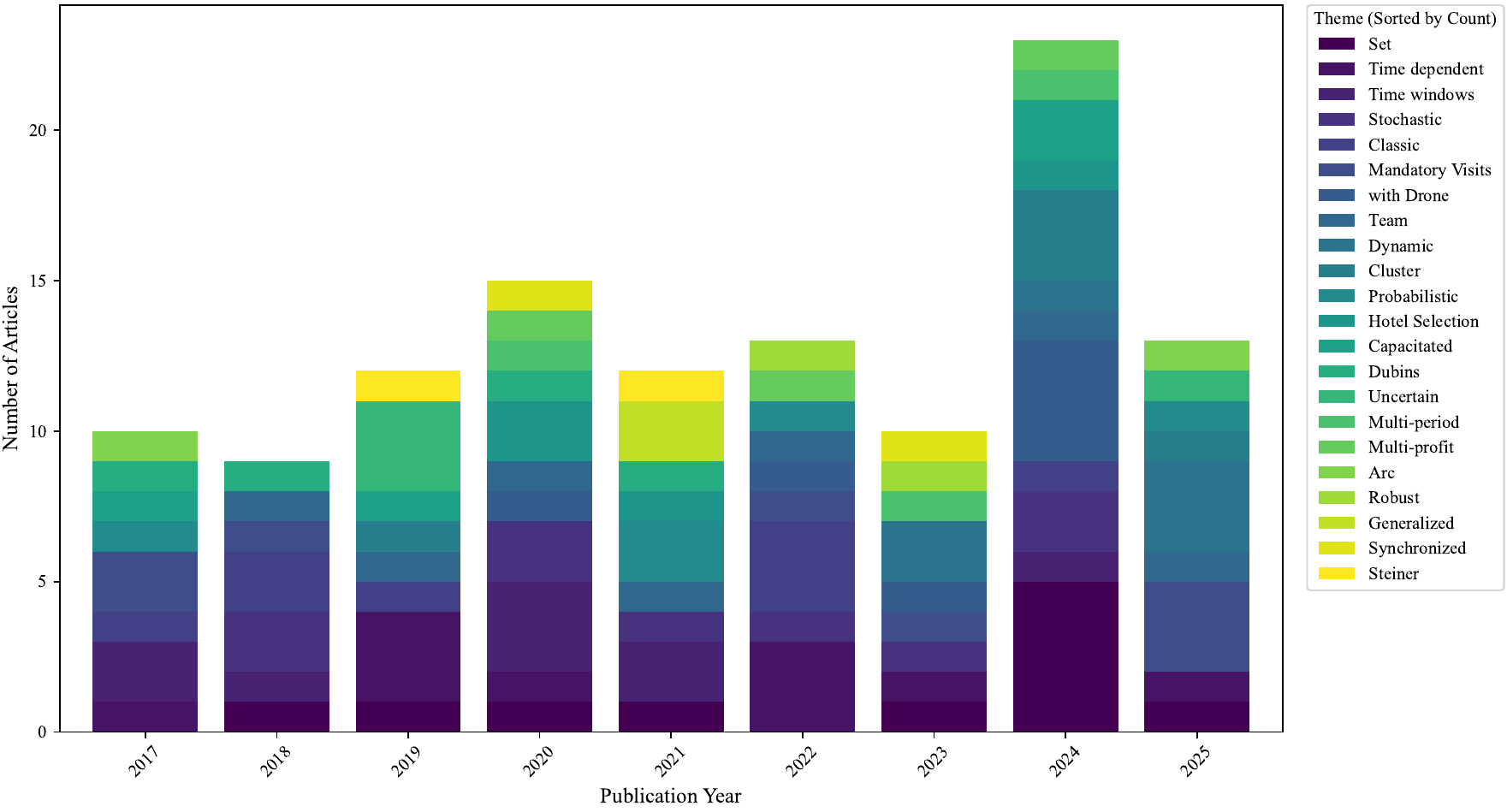}}
	\subfigure[Cumulative trend of articles by theme.]{
		\label{fig:operatorc}
		\includegraphics[width=0.45\linewidth]{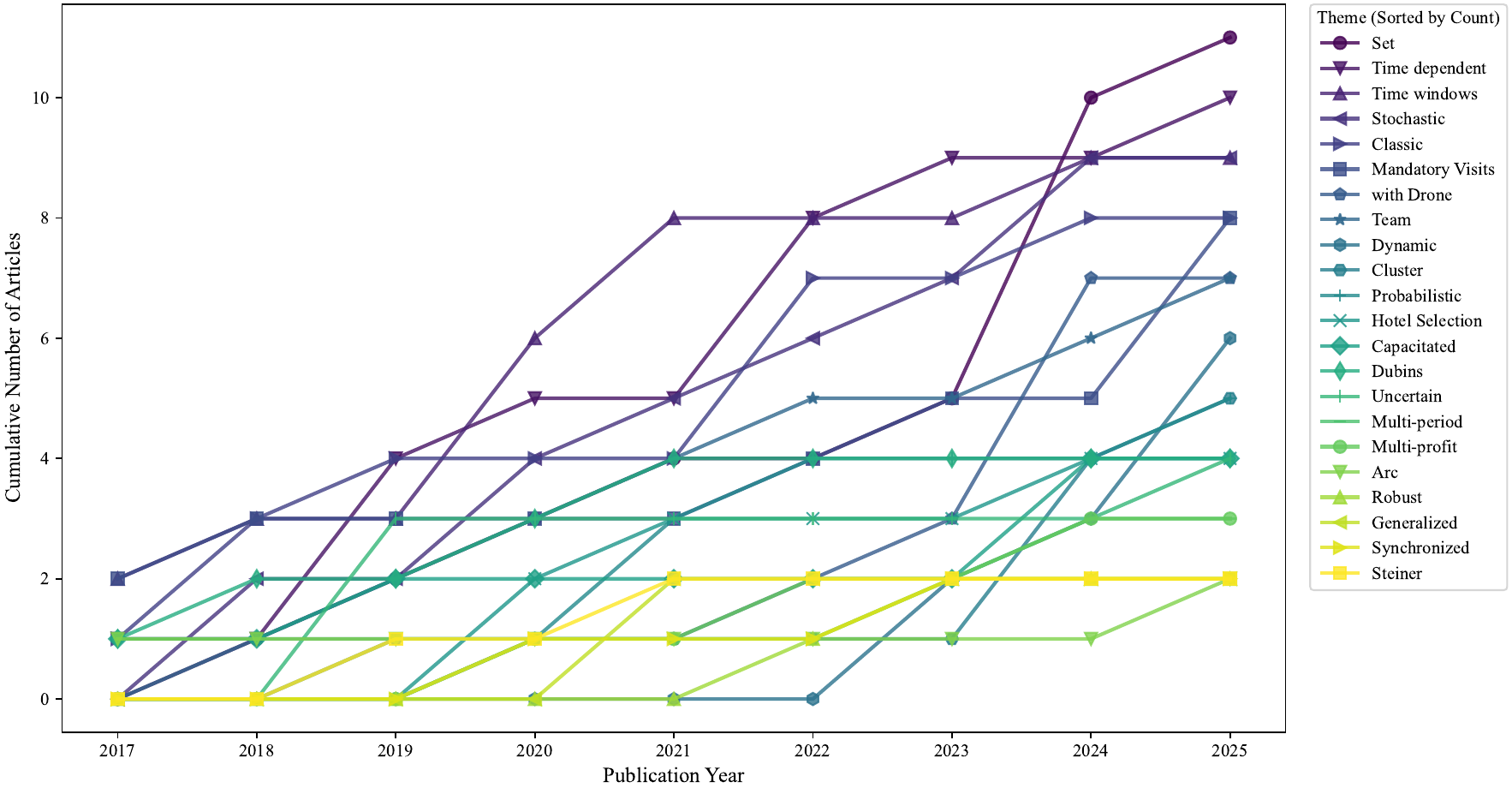}}
	
	\caption{\quad Trends and distribution of the Orienteering Problem research themes (2017–2025).}
	\label{fig:fig2}
\end{figure*}

Research from 2017 to 2025 highlights the dominance of established models like Time-Dependent (TD-OP) and Stochastic OPs (StchOP) for real-time adaptability, alongside emerging variants such as Dynamic (DyOP) and Set OPs (SOP) (Fig.~\ref{fig:fig2}). Methodologically, traditional exact and metaheuristic approaches (e.g., LNS) are increasingly complemented by hybrid strategies combining optimization with Machine Learning. Notably, Reinforcement Learning (RL) and neural-enhanced evolutionary algorithms have become key drivers for improving scalability in dynamic settings.

\subsection{Structural implications of the component-based taxonomy}\label{sec:taxonomy}
Each component in Fig. \ref{fig_fig3} transforms the underlying optimization structure of the Orienteering Problem (OP) in a different way, with consequences for complexity, polyhedral properties and, ultimately, for the design of exact and approximate algorithms.
We view the classic OP as a prize-collecting TSP with a knapsack-type budget, and briefly indicate how each component modifies this baseline.

\textbf{Time-based components.}
Time-based extensions (time windows, time-dependent) can be interpreted as embedding the OP in a time-expanded resource network. Feasibility depends on both the visit sequence and arrival times, and the pseudo-polynomial DP structures of simpler knapsack-like formulations may no longer be directly exploitable. From a polyhedral viewpoint, classical routing inequalities remain valid but are typically combined with time–resource cuts; exact methods therefore tend to rely more on labeling algorithms and branch-and-(cut-and-)price than on compact arc–flow formulations, although the latter are still possible in some settings.

\textbf{Path-based components.}
Path-based extensions (arc OP, Dubins OP) shift the focus from selecting vertices to selecting arcs or motion primitives. The decision space grows on the arc or path level, and geometric constraints (e.g., heading, turning radius) require discretized configuration spaces. The relevant polyhedra are arc- and path-routing polytopes with additional transition constraints, so exact methods naturally adopt path-based set-partitioning and column-generation formulations rather than purely node-based models.

\textbf{Node-based components.}
Node-based extensions (capacity, mandatory vertices) keep the routing topology but enrich vertex constraints. They can be viewed as the intersection of an OP routing polytope with knapsack- and assignment-type polytopes. Classical TSP/OP cuts remain useful, yet effective formulations typically also exploit cover and assignment inequalities; the main structural change lies in tighter local feasibility conditions rather than in a fundamentally different global connectivity pattern.

\textbf{Structure-based components.}
Structure-based extensions (set/clustered OP) modify the combinatorial structure of the vertex set itself. Decisions operate on clusters and on multiple interacting routes, which is closely related to generalized routing and multi-route assignment problems. The resulting polyhedron can be interpreted as combining generalized routing facets with cluster and partitioning facets, which motivates route- and cluster-based set-partitioning models and branch-and-price algorithms in many cases, rather than relying only on simple single-route formulations.

\textbf{Information-based components.}
Information-based extensions (stochastic, dynamic, robust OP) explicitly model uncertainty and information acquisition. They embed OP routes in scenario-based stochastic programs, robust counterparts, or Markov decision processes, where solutions may be formulated as policies rather than fixed routes. The associated extended formulations (scenario expansions, non-anticipativity, uncertainty sets, state transitions) can quickly become very large for pure MILP approaches; as a result, large-scale instances are often addressed with approximate dynamic programming, sampling-based stochastic programming and reinforcement learning techniques.

In summary, the five components correspond to distinct structural transformations of the baseline OP, such as time expansion, arc/path routing, intersection with knapsack-type structures, generalized and clustered routing, and embedding in stochastic or dynamic decision frameworks. This perspective helps explain why different OP variants tend to align with different formulations, cuts and algorithmic paradigms.

\section{Classical models and benchmark variants of the orienteering problem}
This section reviews recent studies on the classic OP and widely studied variants, including the TOP and (T)OPTW.
We focus particularly on works not covered in previous surveys.
The TOP and (T)OPTW are emphasized due to their widespread adoption and integration potential.
Table \ref{tab:Overview of widely studied model} summarizes the models discussed.

As this review focuses on recent developments, we provide only a brief introduction to the basic mathematical model.
For detailed discussions, readers are referred to \cite{vansteenwegen2011}.
Section \ref{sec:3} introduces additional OP variants.

\begin{table}[!tbp]
	\centering
	\caption{Overview of widely studied models and their structural characteristics.}
	\scriptsize	
	\setlength{\tabcolsep}{3pt}	
	\renewcommand{\arraystretch}{1.2}	
	\begin{tabular}{m{2.2cm}<{\raggedright} m{2.0cm}<{\raggedright} m{4.0cm}<{\raggedright} m{4.0cm}<{\raggedright}}	
		\toprule
		\textbf{Specific Problem} & \textbf{Taxonomic Component} & \textbf{Modeling Focus \& Structural Change} & \textbf{Representative Application} \\
		\midrule
		Orienteering Problem (\ref{2.1}) & Baseline Model & Single vehicle; Knapsack-type budget constraint; Prize-collecting objective. & Tourist trip design. \\
		\midrule % 加一条淡线区分，或者用 \addlinespace
		Team Orienteering Problem (\ref{2.2}) & Path-Based Extensions & Multiple vehicles ($m > 1$); Set packing structure; Joint reward maximization. & Fleet management; Multi-agent surveillance; Disaster relief. \\
		\midrule
		OP with Time Windows (\ref{2.3}) & Time-Based Extensions & Visit start time constraints $[O_i, C_i]$; Waiting time permitted; Sequencing dependent feasibility. & Urban logistics; Technician routing; Satellite scheduling. \\
		\bottomrule
	\end{tabular}%
	\label{tab:Overview of widely studied model}%
\end{table}%

\subsection{Orienteering Problem}\label{2.1}
The Orienteering Problem is a combinatorial optimization problem that determines an optimal path through a set of locations.
The objective is to maximize collected rewards subject to resource constraints (e.g., time or distance).
The OP finds applications in route planning, tourism, and logistics, optimizing traversal under resource limits.

\begin{table}[htbp]
	\footnotesize
	\centering
	\caption{Summary of Orienteering Problem algorithms.}
	\scriptsize % 1. 使用小号字体 (如果还觉得大，可以改成 \tiny)
	\setlength{\tabcolsep}{1.5pt} % 2. 极大地减小列间距 (默认是 6pt)
	\renewcommand{\arraystretch}{0.85}
	\begin{tabular}{m{2.5cm}<{\centering} m{2cm}<{\centering} m{5cm}<{\centering} m{3cm}<{\centering}}
		\toprule
		\textbf{Reference} & \textbf{Problem} & \textbf{Algorithm} & \textbf{Performance} \\
		\midrule
		\cite{kobeaga2024}	 & OP & Enhanced branch-and-cut with new bounds & 18 new optima, 76 new BKS, 85 new UBs. \\
		\cite{gottlieb2022} & OP & PTAS for high-dimensional problems & Reduced time complexity to $n^{O(\frac{1}{\delta})}$. \\
		\cite{santini2019} & OP & Adaptive Large Neighborhood Search (ALNS) heuristic & 27 new best-known solutions achieved. \\
		\cite{ostrowski2017} & OP & Evolution-Inspired Local Improvement Algorithm (EILIA) & Solved instances with up to 908 nodes. \\
		\cite{tregel2021} & OP & Randomized Greedy + VNS & Real-time route calculation \\		
		\cite{elzein2022} & OP & Cluster-based decomposition heuristic & Solved problems with over 1,000 nodes. \\
		\cite{kobeaga2018} & OP & Evolutionary algorithm with edge recombination & Tested on 344 instances; strong on large-scale. \\
		\cite{sun2022} & OP & ML-ACO: Machine Learning-augmented ACO & Improved efficiency for 1,000+ nodes. \\
		\bottomrule
	\end{tabular}
	\label{tab:opra}
\end{table}

Formally, the OP is defined on a graph.
Let $G=(V,E)$ be a weighted, directed graph, where $V$ is the set of nodes (locations).
$E$ is the set of edges (paths). Each edge $(i,j) \in E$ has a cost $c_{ij}$ (or $t_{ij}$), representing distance or travel time.
Each node $v_i \in V$ has a reward $S_i$.

The objective is to determine a path starting at node $s$, visiting a subset of nodes, and ending at a destination node, maximizing the total collected reward.
The solution must satisfy a total path cost constraint, typically a time or distance limit.

Following the notation in \cite{vansteenwegen2011}, let $X_{ij} = 1$ if node $i$ is followed by node $j$, and 0 otherwise.
The variable $U_i$ represents the position of the visited node $i$ in the path, serving to enforce subtour elimination constraints.

Therefore, the OP can be described by an integer programming model as:

\vspace{-2.0em}
\begin{align}
	\label{objfun}
	& \max \sum_{i = 2}^{|V| - 1} \sum_{j = 2}^{|V|}S_{i}X_{ij} \\
	\label{opc1}
	& \sum_{j = 2}^{|V|}X_{1j} = \sum_{i = 1}^{|V| - 1}X_{i|V|} = 1 \\
	\label{opc2}
	& \sum_{i = 1}^{|V| - 1} X_{ik} = \sum_{j = 2}^{|V|} X_{kj} \leq 1, \quad \forall k = 2, \dots, (|V| - 1) \\
	\label{opc3}
	& \sum_{i = 1}^{|V| - 1} \sum_{j = 2}^{|V|} t_{ij}X_{ij} \leq T_{max} \\
	\label{opc4}
	& 2 \leq U_i \leq |V|, \quad \forall i = 2, \dots, |V| \\
	\label{opc5}
	& U_i - U_j + 1 \leq (|V| - 1)(1 - X_{ij}), \quad \forall i = 2, \dots, |V|
\end{align}

The objective function maximizes the total collected score.
The interpretation of the constraints is as follows:
Equation (\ref{opc1}) ensures that the solution path starts at node 1 and ends at node \(|V|\). 
Equation (\ref{opc2}) enforces flow conservation and ensures intermediate nodes are visited at most once.
Specifically, it ensures that for all nodes except node 1 and node \(|V|\), the in-degree equals the out-degree, and that each node is visited at most once. 
Equation (\ref{opc3}) restricts the total path cost to a predefined limit $T_{max}$.
Equations (\ref{opc4}) and (\ref{opc5}) represent the Miller-Tucker-Zemlin (MTZ) constraints, used to eliminate subtours.

Table \ref{tab:opra} summarizes recent studies on the classic OP.

\subsubsection{Benchmark instances}
The standard benchmarks rely on TSPLib-based instances generated by \cite{fischetti1998} (Generations 1--3). \cite{kobeaga2024} introduced a fourth generation with varying distance constraints to specifically challenge Branch-and-Cut (B\&C) algorithms. These datasets, along with best-known solutions from EA4OP, are publicly available on GitHub.

\subsubsection{State-of-the-art results}

The development of solution methods for the OP reflects a persistent tension between the scalability of heuristics and the precision of exact algorithms. Early research largely deemed exact methods unsuitable for large-scale instances, but recent work challenges this view via advanced decomposition and pricing techniques. In parallel, heuristic approaches have evolved from simple constructive rules to sophisticated metaheuristics that exploit infeasibility, spatial structure, and, more recently, Machine Learning.

\cite{kobeaga2024} proposed a Revisited Branch-and-Cut (RB\&C) incorporating shrinking strategies and efficient variable pricing. This approach solved instances with up to 7,397 nodes, yielding 18 new optimal solutions and 76 Best Known Solutions (BKS), demonstrating that advanced pricing can unlock large-scale exact optimization.

On the heuristic side, Exploiting infeasible regions has proved effective. \cite{ostrowski2017} (EILIA) and \cite{kobeaga2018} (EA4OP) both maintain infeasible intermediates to escape local optima. EA4OP, which couples edge-recombination crossover with Lin–Kernighan refinement, demonstrates superior robustness on large-scale instances compared to prior state-of-the-art.

\cite{santini2019} introduced an Adaptive Large Neighborhood Search (ALNS) using clustering-based operators to exploit topological structure, achieving 27 new BKS. A comparison reveals complementarity: EA4OP \citep{kobeaga2018} excels in speed, while ALNS \citep{santini2019} often yields higher quality on disjoint hard instances.

\cite{sun2022} integrated Machine Learning into Ant Colony Optimization (ML-ACO). By training classifiers on small graphs to guide pheromone updates, they demonstrated that learned patterns generalize effectively to large synthetic and real-world networks.

Tables \ref{tab:op1} and \ref{tab:op2} in \ref{append1} detail the new Best Known Solutions found in recent research.

\subsection{Team Orienteering Problem}\label{2.2}
The Team Orienteering Problem (TOP) extends the classic OP within combinatorial optimization.
It involves determining paths for a set of teams. 
The goal is to maximize total rewards while satisfying resource constraints for each path.
Applications include multi-vehicle route planning, multi-agent systems, and collaborative logistics. 
In these scenarios, efficient resource distribution is critical.

\begin{table}[htbp]\label{tab:top_summary}
	\footnotesize
	\centering
	\caption{Summary of methods for solving the Team Orienteering Problem.}
	
	\scriptsize % 1. 使用小号字体 (如果还觉得大，可以改成 \tiny)
	\setlength{\tabcolsep}{1.5pt} % 2. 极大地减小列间距 (默认是 6pt)
	\renewcommand{\arraystretch}{0.85}
	\begin{tabular}{m{2.5cm}<{\centering} m{2cm}<{\centering} m{3cm}<{\centering} m{5cm}<{\centering}}			
		\toprule
		\textbf{Reference} & \textbf{Problem} & \textbf{Algorithm} & \textbf{Performance} \\
		\midrule
		\cite{bianchessi2018} & TOP & B\&C & Solved 327 (small) instances, 24 new \\
		\cite{hammami2020} & TOP & HALNS & Best solutions for 720 instances (1 new) \\			
		\cite{tsakirakis2019} & TOP & SHHS & Best-known solutions for 276/328 instances \\
		\cite{panadero2021} & TOP & Biased Randomization & Real-time solutions ($<1$ second) \\
		\cite{peyman2024} & TOP & Sim-Learnheuristic & High-quality solutions rapidly \\
		\cite{chaigneau2025} & VLS-TOP & LNS + Clustering & 3 new BKS, Scalable to 5000+ nodes \\
		\bottomrule
	\end{tabular}
	\label{tab:topsum}
\end{table}

Formally, the TOP is defined on a graph.
Let $G=(V,E)$ be a weighted, directed graph, where:

$V$ is the set of nodes (locations).

$E$ is the set of edges. Each edge $(i,j) \in E$ has a cost $c_{ij}$ (or $t_{ij}$), representing distance or travel time.

Each node $v_i \in V$ has a reward $S_i$.

The objective is to determine $P$ paths, each starting at node $s$ and ending at a destination node. 
The solution must maximize the total collected reward.
Each path must satisfy a cost constraint (e.g., a time or distance limit).

Let $X_{ijp} = 1$ if path $p$ traverses edge $(i, j)$, and 0 otherwise. 
Let $Y_{ip} = 1$ if path $p$ visits node $i$. 
Variable $U_{ip}$ represents the position of node $i$ in path $p$, used for subtour elimination.

The Integer Programming (IP) model is formulated as follows:

\vspace{-2.0em}
\begin{align}
	\label{objfun_top}
	& \max \sum_{p = 1}^{P} \sum_{i = 2}^{|V| - 1}S_{i}Y_{ip} \\
	\label{topc1}
	& \sum_{p = 1}^{P}\sum_{j = 2}^{|V|}X_{1jp} = \sum_{p = 1}^{P}\sum_{i = 1}^{|V| - 1}X_{i|V|p} = 1 \\
	\label{topc2}
	& \sum_{p = 1}^{P}Y_{kp}\leq 1, \quad \forall k = 2, \dots, (|V| - 1) \\
	\label{topc3}
	& \sum_{i = 1}^{|V| - 1} X_{ikp} = \sum_{j = 2}^{|V|} X_{kjp} = Y_{kp}, \quad \forall k = 2, \dots, (|V| - 1), \forall p = 1,\dots,P \\
	\label{topc4}
	& \sum_{i = 1}^{|V| - 1} \sum_{j = 2}^{|V|} t_{ij}X_{ijp} \leq T_{max}, \quad \forall p = 1,\dots,P \\
	\label{topc5}
	& 2 \leq U_{ip} \leq |V|, \quad \forall i = 2, \dots, |V|, \forall p = 1,\dots,P \\
	\label{topc6}
	& U_{ip} - U_{jp} + 1 \leq (|V| - 1)(1 - X_{ijp}), \quad \forall i = 2, \dots, |V|, \forall p = 1,\dots,P
\end{align}

The objective function maximizes the total score. 
The constraints are defined below:
Equation (\ref{topc1}) ensures that all $P$ paths in the solution start at node 1 and end at node \(|V|\). 
Equation (\ref{topc2}) limits each intermediate node, except for nodes \(1\) and \(|V|\), to be visited by at most one path.
Equation (\ref{topc3}) enforces flow conservation.
Equation (\ref{topc4}) restricts the cost of each path to the limit $T_{max}$.
Equations (\ref{topc5}) and (\ref{topc6}) implement MTZ constraints for subtour elimination.

Table \ref{tab:topsum} summarizes recent TOP research.

\subsubsection{Benchmark instances}
Standard studies utilize the \cite{chao1996} benchmarks (up to 401 nodes). Recently, to address industrial needs, \cite{chaigneau2025} introduced large-scale instances (1,001–5,395 customers) adapted from \cite{kobeaga2018}, explicitly testing algorithmic scalability limits.

\subsubsection{State-of-the-art results}

\cite{bianchessi2018} challenged the dominance of column generation by proposing a compact two-index formulation with dynamic connectivity cuts. This polynomial-sized model closed 24 previously open cases among small-scale benchmarks. However, the cost of exact separations restricts its applicability to graphs under $\sim$100 nodes. To balance quality and scale, \cite{hammami2020} proposed the Hybrid Adaptive Large Neighborhood Search (HALNS). Uniquely, it embeds an exact Set Packing model to optimally recombine routes, achieving 100\% of Best Known Solutions (BKS) on standard sets. While superior in precision, the reliance on a MIP solver creates a bottleneck as route pools expand. Addressing industrial scalability, \cite{chaigneau2025} shifted to a ``cluster-first, route-second'' decomposition. By pruning the search space via balanced k-medoids, this LNS approach handles 7,000-node instances within minutes and discovered 3 new BKS. Although it yields slightly higher gaps on medium instances than the MIP-hybrid HALNS, it dominates in ultra-large-scale regimes where global optimization is intractable.

Table \ref{tab:nbks_top} (\ref{append1}) details the new BKS for TOP instances.

\subsection{(Team) Orienteering Problem with Time Windows}\label{2.3}
The OP with Time Windows (OPTW) and its extension, the Team Orienteering Problem with Time Windows (TOPTW), integrate time window constraints into the OP framework.
In the OPTW, each node is assigned a time window $[O_i, C_i]$, defining when service must start.
Early arrival necessitates waiting, while late arrival makes the visit infeasible.
\cite{kantor1992} first investigated the single-path OPTW ($m = 1$).
Extending the problem to multiple paths ($m > 1$) yields the TOPTW.

\begin{table}[htbp]
	\centering
	\caption{Summary of methods for solving the (T)OPTW.}	
	\label{tab:toptw_summary}	
	\scriptsize % 1. 使用小号字体 (如果还觉得大，可以改成 \tiny)
	\setlength{\tabcolsep}{1.5pt} % 2. 极大地减小列间距 (默认是 6pt)
	\renewcommand{\arraystretch}{1}
	\begin{tabular}{m{2.5cm}<{\centering} m{2cm}<{\centering} m{3cm}<{\centering} m{5cm}<{\centering}}
		\hline
		\textbf{Reference} & \textbf{Problem} & \textbf{Algorithm} & \textbf{Performance} \\ \hline
		\cite{gedik2017} & TOPTW & Constraint Programming & 119 new best-known solutions \\
		\cite{gunawan2017} & TOPTW & SA \& ILS & Improved 50 (SA) / 37 (ILS) benchmarks \\
		\cite{gunawan2018} & TOPTW & ADOPT Framework & 11 new best solutions \\
		\cite{karabulut2020} & TOPTW & Evolutionary Strategy & 7 new best solutions \\
		\cite{amarouche2020} & TOPTW & ILS + Adaptive Memory & 57 new high-quality solutions \\	
		\cite{gama2021} & OPTW & Pointer Network + RL & Better than heuristics in benchmarks \\
		\cite{moosavi2022} & TOPTW & MOGA & Superior to NSGA-II on large-scale \\
		\cite{tran2024} & OPTW & Route Recombination & Efficient post-optimization \\
		\hline			
	\end{tabular}
\end{table}

Using the previous notation, the TOPTW is formulated as a Mixed-Integer Linear Programming (MILP) model:
Let $X_{ijp} = 1$ if path $p$ traverses edge $(i, j)$, and 0 otherwise. 
Let $Y_{ip} = 1$ if path $p$ visits node $i$.
Variable $W_{ip}$ represents the service start time at node $i$ in path $p$.

\vspace{-2.0em}
\begin{align}
	\label{objfun_toptw}
	& \max \sum_{p = 1}^{P} \sum_{i = 2}^{|V| - 1}S_{i}Y_{ip} \\
	\label{toptwc1}
	& \sum_{p = 1}^{P}\sum_{j = 2}^{|V|}X_{1jp} = \sum_{p = 1}^{P}\sum_{i = 1}^{|V| - 1}X_{i|V|p} = P \\
	\label{toptwc2}
	& \sum_{p = 1}^{P}Y_{kp}\leq 1, \quad \forall k = 2, \dots, (|V| - 1) \\
	\label{toptwc3}
	& \sum_{i = 1}^{|V| - 1} X_{ikp} = \sum_{j = 2}^{|V|} X_{kjp} = Y_{kp}, \quad \forall k = 2, \dots, (|V| - 1), \forall p = 1,\dots,P \\
	\label{toptwc4}
	& W_{ip} + t_{ij} - W_{jp} \leq M(1 - X_{ijp}), \quad \forall i,j = 1, \dots, |V|, \forall p = 1,\dots,P \\
	\label{toptwc5}
	& W_{|V|p} - W_{1p} \leq T_{\max}, \quad \forall p = 1,\dots,P \\
	\label{toptwc6}
	& O_i \leq W_{ip}, \quad \forall i = 1, \dots, |V|, \forall p = 1,\dots,P \\
	\label{toptwc7}
	& W_{ip} \leq C_i, \quad \forall i = 1, \dots, |V|, \forall p = 1,\dots,P
\end{align}
\vspace{-2.0em}

The objective function (\ref{objfun_toptw}) maximizes the total collected score.
Equation (\ref{toptwc1}) ensures that all $P$ paths in the solution start at node 1 and end at node \(|V|\). 
Equation (\ref{toptwc2}) limits each intermediate node, except for nodes \(1\) and \(|V|\), to be visited by at most one path.
Equation (\ref{toptwc3}) enforces flow conservation.
Equation (\ref{toptwc4}) ensures timeline consistency between consecutive nodes.
Equation (\ref{toptwc5}) restricts the duration of each path to the limit $T_{max}$.
Excluding waiting time, the constraint simplifies to:
\vspace{-1.0em}
\begin{equation*}
	\sum_{i = 1}^{|V| - 1} \sum_{j = 2}^{|V|} t_{ij}X_{ijp} \leq T_{max}; \forall p = 1,\dots,P.
\end{equation*}
\vspace{-1.0em}

Equations (\ref{toptwc6}) and (\ref{toptwc7}) enforce time window constraints.

Setting $P=1$ reduces this model to the single-path OPTW.

\subsubsection{Benchmark instances}
Reviewed studies utilize benchmark datasets based on modified existing instances to test OP solution approaches.
For instance, \cite{gedik2017} used 304 TOPTW instances derived from \cite{solomon1987} and \cite{cordeau1997}.
These instances were adapted to include constraints reflecting real-world routing scenarios.

\subsubsection{State-of-the-art results}
\cite{gedik2017} established a rigorous baseline using CP with tailored branching on interval variables, discovering 119 Best Known Solutions (BKS) and certifying optimality for two open instances. However, the scalability limits of exact proofs on large benchmarks necessitated a shift toward heuristic adaptivity.

Initial metaheuristics, such as the hybrid SA-ILS by \cite{gunawan2017}, relied on offline parameter tuning to achieve 50 BKS, but suffered from static calibration. To address search dynamics, \cite{gunawan2018} introduced ADOPT, utilizing online learning automata to update operator probabilities in real-time, while \cite{karabulut2020} employed Evolution Strategies with self-adaptive mutation. These dynamic frameworks systematically outperformed static designs by evolving perturbation strengths alongside the search process.

A distinct structural shift was introduced by \cite{amarouche2020}. Departing from direct route construction (common in CP and SAILS), they adopted a "giant tour" decomposition. An optimal splitting procedure based on an $O(n^2)$ dynamic program decouples customer sequencing from capacity and time-window constraints. This indirect search achieved a 0.30\% optimality gap and 57 BKS, suggesting that decoupling feasibility checks from the search space offers a superior exploration–exploitation balance for highly constrained instances compared to direct construction methods.

Tables \ref{tab:toptw1} and \ref{tab:toptw2} in \ref{append1} detail the new BKS for TOPTW instances.

\subsection{Discussion}
Foundational OP research (OP, TOP, TOPTW) is increasingly blurring the traditional exact-heuristic dichotomy. Metaheuristics now integrate exact components, exemplified by HALNS \citep{hammami2020}, which embeds a MIP-based set packing model within large neighborhood search to optimize route selection. Conversely, exact algorithms like the Revisited Branch-and-Cut \citep{kobeaga2024} have adopted heuristic-inspired shrinking and pricing strategies, demonstrating that engineered exact schemes can still dominate massive deterministic instances.

In parallel, solution representations have evolved to bridge the gap between academic benchmarks and industrial scale. By shifting search spaces from direct edge selection to permutations (giant tours) \citep{amarouche2020} or clusters \citep{chaigneau2025}, recent methods structurally decouple feasibility checks from resource constraints. This decomposition allows solvers to handle thousands of customers—a scale previously unreachable by direct construction methods. These hybrid, decomposition-based frameworks provide the computational backbone for the complex stochastic and dynamic extensions discussed next.

\section{Advances in canonical extensions of the orienteering problem}\label{sec:3}

\begin{table}[!th]
	\centering
	\caption{Overview of canonical extended models.}
	\scriptsize	
	\setlength{\tabcolsep}{3pt}	
	\renewcommand{\arraystretch}{1.2}	
	\begin{tabular}{m{2.2cm}<{\raggedright} m{2.0cm}<{\raggedright} m{4.0cm}<{\raggedright} m{4.0cm}<{\raggedright}}	
		\toprule
		\textbf{Specific Problem} & \textbf{Taxonomic Component} & \textbf{Modeling Focus \& Structural Change} & \textbf{Representative Application} \\
		\midrule
		
		Stochastic Orienteering Problem (\ref{3.1}) & Information-Based Extensions & Modeling parameters as random variables; optimization targets expected value or satisfies chance constraints. & Robotics navigation, search and rescue operations, etc. \\
		
		Time-Dependent Orienteering Problem (\ref{3.2}) & Time-Based Extensions & Cost/reward defined as functions of arrival time ($c_{ij}(t)$); introduces non-linear temporal consistency. & Traffic-sensitive urban logistics, time-specific tourist route planning, etc. \\
		
		Capacitated Orienteering Problem (\ref{3.4}) & Node-Based Extensions & Adding demand accumulation constraints; routes bounded by both travel budget and vehicle capacity. & Disaster relief logistics involving constrained resources, etc. \\
		
		Orienteering Problem with Mandatory Visits (\ref{3.3}) & Node-Based Extensions & Hybridizing OP with TSP/Steiner structures; imposes hard feasibility constraints on specific node subsets. & Priority customer services, critical infrastructure inspection, etc. \\	
		
		Clustered Orienteering Problem (\ref{3.5}) & Structure-based extensions & Profit associated with node subsets (clusters); enforcing atomic visitation or specific intra-cluster rules. & Group-buying logistics, etc. \\
		
		\bottomrule
	\end{tabular}%
	\label{tab:Overview of established Extended model}%
\end{table}%

This section reviews recent developments in canonical extensions of the orienteering problem, representing well-established problem families previously covered in surveys such as \cite{vansteenwegen2011} and \cite{gunawan2016}. Our focus here is on methodological advances and new applications reported in the recent period.

Although Section~\ref{sec:taxonomy} introduces a component-based taxonomy, we organize this section by established problem names (e.g., Stochastic OP). This reflects the maturity of these families, where research now prioritizes algorithmic refinement over model formulation. Unlike the emerging models in Section~\ref{sec:4}—where structural categorization is essential to clarify novel constraints—canonical variants are best navigated by traditional nomenclature. Table~\ref{tab:Overview of established Extended model} maps these variants to our taxonomic framework.

\subsection{Stochastic Orienteering Problem}\label{3.1}

Building on early surveys \cite{gunawan2016}, we review recent advances in the Stochastic Orienteering Problem (StchOP). StchOP extends the deterministic model by treating node profits, travel, and service times as stochastic variables, requiring the optimization of expected performance or probabilistic objectives \cite{ilhan2008}.

\subsubsection{Solution approaches}
A central challenge in StchOP is the “myopia” of traditional re-optimization strategies, which react to realized uncertainty without anticipating future stochastic outcomes. \cite{bian2018} address this issue with a Simulation-Aided Multiple Plan Approach (SMPA) that generates a pool of diverse routes during idle periods and selects among them as uncertainty unfolds. By explicitly embedding a “myopia-prevention” mechanism, SMPA yields markedly more robust policies than standard single-plan re-optimization. Complementing this operational focus, \cite{dolinskaya2018} incorporate adaptivity directly into strategic planning by combining Dynamic Programming with a Variable Neighborhood Search. Their results show that modeling recourse actions already at the node-selection stage systematically outperforms static planning, highlighting the value of anticipative design rather than purely reactive control.

Beyond expected-profit maximization, recent work has explored richer stochastic objectives and service-related constraints. \cite{panadero2022} study the Stochastic Team OP with Position-Dependent Rewards (StchTOP-PDR), where travel-time variability directly affects reward values. They show that a simheuristic framework—metaheuristics coupled with simulation—scales more effectively than exact methods for such highly nonlinear dependencies, indicating a pragmatic trade-off between optimality and tractability. In home-delivery applications, \cite{Song2020} formulate a two-stage stochastic model that enforces driver–customer consistency, treating “trust” as an explicit objective alongside profit. Their multiple-scenario approach captures the structural rigidity induced by consistency constraints and stands in contrast to purely dynamic routing schemes that privilege flexibility. Finally, risk-sensitive variants have adopted chance constraints to control failure probabilities, particularly in robotic settings where budget exhaustion is critical. \cite{thayer2021} propose an adaptive policy based on Constrained Markov Decision Processes (CMDP), but \cite{carpin2024} show that CMDP formulations suffer from state-space explosion due to time discretization. Their online Monte Carlo Tree Search (MCTS) algorithm avoids discretization and thus provides a more scalable approach to continuous-state, chance-constrained StchOP.

\subsubsection{Discussion}
Reflecting the evolution of information-based extensions, work on the StchOP has shifted from static, expectation-based planning to designs that treat uncertainty as a first-class modeling concern. One line of research embeds simulation directly into the optimization loop—through SMPA and simheuristics—using sampling to assess the downstream impact of routing decisions and thereby mitigating the myopia of single-plan re-optimization. A complementary line emphasizes explicit adaptivity and risk control, from node-selection schemes that anticipate recourse to chance-constrained CMDP formulations and online MCTS policies. Together, these approaches move the focus from optimizing a nominal plan to managing operational risk over time, enabling objectives that capture service quality and trust while remaining solvable on realistic instances.

\subsection{Time-Dependent Orienteering Problem}\label{3.2}
The Time-Dependent Orienteering Problem (TD-OP) extends the classic OP by incorporating temporal variations in travel times, service durations, and rewards.
Unlike static OP frameworks, the TD-OP captures dynamic conditions like traffic congestion, increasing both complexity and practical relevance.
Foundational work by \cite{fomin2002} introduced approximation algorithms for time-dependent scenarios, establishing a basis for subsequent research.
Recent advancements focus on specialized algorithms and applications in transportation, tourism, and satellite scheduling, emphasizing the role of temporal dynamics.

\subsubsection{Solution approaches}
The Time-Dependent Orienteering Problem (TD-OP) extends static routing models by recognizing that both travel conditions and rewards may vary over time. Early work mainly focused on congestion: \cite{verbeeck2017} studied the Time-Dependent Orienteering Problem with Time Windows (TD-OPTW) and proposed a fast Ant Colony System capable of solving realistic road-network instances within seconds, setting an early benchmark for congestion-aware routing. Subsequent studies argued that time dependence should also apply to rewards. \cite{yu2019a} introduced the Orienteering Problem with Service Time Dependent Profits (OPSTP), where profits grow with service duration, and solved instances of up to 200 vertices via a two-phase matheuristic combining Tabu Search and an exact scheduling routine. In a related direction, \cite{yu2019b} considered arrival-time-dependent scores and proposed a Hybrid Artificial Bee Colony algorithm with a Simulated Annealing acceptance rule that outperforms standard metaheuristics on medium and large instances.

A more integrated view of time dependence is offered by \cite{khodadadian2022}, who combine time-dependent travel times with service-time-dependent profits in TDOPTW-STP and solve it via Variable Neighborhood Search. Their results indicate that treating congestion and flexible service durations in isolation underestimates potential efficiency gains in congested networks. Increasingly complex profit structures further challenge traditional heuristic practice. \cite{yu2021} study a variant with non-concave profits depending jointly on arrival time and service duration, and, in contrast to earlier heuristic-only approaches, develop an exact Benders Branch-and-Cut algorithm alongside a hybrid heuristic (ILS-MCS) that handles instances with up to 100 customers. A similar trend toward exact methods is evident in satellite scheduling: \cite{peng2019} first address agile satellite scheduling using a heuristic based on bidirectional dynamic programming and Iterated Local Search, then move to an adaptive-directional dynamic programming approach for single-orbit problems \cite{peng2020}. Most recently, \cite{peng2025} propose a Branch-and-Cut-and-Price algorithm with subset-row inequalities for the Time-Dependent Team Orienteering Problem with Variable Time Windows (TD-TOPVTW), achieving more than a sixfold speedup over state-of-the-art exact algorithms on standard TOPTW benchmarks.

Beyond congestion and simple time-dependent profits, several structural extensions highlight the breadth of TD-OP modeling. \cite{dasdemir2022} study UAV routing with time-dependent prizes, preserving static travel times but offering multiple trajectory options between nodes to trade off risk and duration, solved by a hybrid heuristic–MIP approach. Addressing the limitations of deterministic time windows, \cite{avraham2023} adopt a data-driven model with soft time windows and stochastic travel and service times, using a Branch-and-Bound algorithm enhanced by local search. Their results underscore that hard time windows may be overly rigid in uncertain environments, and that robustness considerations are essential for realistic TD-OP applications.

\subsubsection{Discussion}
Recent work on Time-Dependent OPs has progressed from congestion-focused models with time-dependent travel times to richer formulations where service durations and rewards also depend on time, showing that treating these components jointly can substantially alter route profitability and feasibility \citep{yu2019a,yu2019b,khodadadian2022,yu2021}. Algorithmically, there is a clear trend from fast heuristics toward matheuristics and exact methods—such as hybrid local-search–MIP schemes, Benders Branch-and-Cut, and Branch-and-Cut-and-Price—that can handle complex temporal structures at moderate scales, while structurally enriched variants for UAV routing and stochastic, data-driven soft time windows emphasize robustness in uncertain environments \citep{dasdemir2022,avraham2023,peng2025}. Overall, TD-OP formulations now serve as a flexible framework for dynamic routing in congested road networks, tourist trip design, UAV operations, and satellite scheduling, balancing model realism with scalable solutions, thereby solidifying the time-based dimension as a critical resource in modern routing.

\subsection{Orienteering Problem with Mandatory Visits}\label{3.4}
The foundational concept originated with \cite{gendreau1998}, who proposed the Undirected Selective Traveling Salesman Problem.
The Orienteering Problem with Mandatory Visits (OPMV) combines mandatory and optional node visits, reflecting applications in logistics, tourism, and disaster management.
While early research focused on basic mandatory constraints, recent advancements have significantly expanded the OPMV's scope.
Contemporary studies tackle larger scales using Machine Learning (ML) \citep{fang2023} and IoT data \citep{li2022}. They also integrate complex features like node incompatibilities \citep{perez2025, guastalla2025} and simultaneous production \citep{wang2025}.
Furthermore, variants like the Steiner Team Orienteering Problem (SteinerTOP) remain significant for their theoretical properties \citep{assuncao2019}.

\subsubsection{Solution approaches}
Unlike the classic OP, where an empty route is technically feasible, finding a valid solution in OPMV variants is already challenging, as the feasibility problem itself can be NP-hard \citep{assuncao2021}. This has led to distinct methodological streams, from advanced heuristics that seek scalability, to exact methods targeting structural complexity, and, more recently, dynamic and data-driven approaches.

Early work relied mainly on metaheuristics to cope with mandatory visits combined with additional constraints. For the OPMVC, \cite{palomo2017} proposed a hybrid GRASP–VNS, while \cite{lin2017} used a Multi-start Simulated Annealing for the TOPTW-MV to improve over both exact solvers and earlier metaheuristics. \cite{lu2018} further refined solution quality with a memetic algorithm that couples Tabu Search with backbone-based crossover. Yet these heuristic frameworks can still struggle to identify feasible regions in highly constrained instances. To address this, \cite{assuncao2021} introduced a matheuristic for the SteinerTOP that combines Large Neighborhood Search with a Feasibility Pump, explicitly exploiting linear relaxations to drive the search toward feasibility when traditional constructive methods fail.

In parallel, exact methods have re-emerged as viable tools for structurally complex OPMV variants where heuristic feasibility is fragile. \cite{assuncao2019} formalized the SteinerTOP via a compact commodity-based formulation and a cutting-plane algorithm effective on medium-sized instances. Extending this line, \cite{perez2025} showed that an undirected Branch-and-Cut for the OPMVC can solve cases previously addressed only by heuristics, while \cite{guastalla2025} proved NP-completeness for TOP-ST-MIN feasibility under mandatory and incompatible nodes and proposed cycle-based cuts that outperform standard Branch-and-Cut implementations. In additive manufacturing, \cite{wang2025} highlight that exact methods become particularly attractive when routing must be synchronized with non-spatial constraints such as production schedules, where the variability of heuristic solutions is less acceptable.

A third stream emphasizes responsiveness and data-driven decision-making in dynamic environments. \cite{li2022} integrated IoT analytics with agile optimization for a dynamic waste collection problem modeled as a Dynamic TOP-MV, showing that re-optimizing routes with real-time information can substantially increase collected reward relative to static baselines. Going further, \cite{fang2023} designed a neural-network-based heuristic for the TOPMV using a Hierarchical Recurrent Graph Convolutional Network within a Large Neighborhood Search, demonstrating that learned operators can outperform hand-crafted heuristics on large-scale instances. These results suggest a promising direction in which learned heuristics not only enhance scalability but may also serve as effective warm starts for exact algorithms in OPMV settings.

\subsubsection{Discussion}
OPMV research has matured along three complementary lines driven by node-based constraints. While metaheuristics (e.g., GRASP–VNS) address general feasibility, the resurgence of exact methods (e.g., branch-and-cut for SteinerTOP) proves that structural node properties can be exploited to certify optimality where heuristics fail. Concurrently, data-driven approaches demonstrate how learned operators sustain quality in dynamic settings. Together, these strands position mandatory visits not merely as constraints, but as structural drivers for algorithmic design.

\subsection{Capacitated Orienteering Problem}\label{3.3}
The Capacitated Orienteering Problem (COP) extends the classic OP by incorporating vehicle capacity constraints, increasing the optimization challenge.
By balancing profit maximization with capacity limitations, the COP addresses applications in logistics, tourism, and disaster management.
Early work by \cite{archetti2009} established the theoretical foundation and proposed initial methods. 
Subsequent research has advanced these models using hybrid heuristics and adaptive algorithms, with a growing focus on computational performance and real-world applications.

\subsubsection{Solution approaches}

Establishing a rigorous foundation, \cite{park2017} presented the first exact algorithm for the CTOP with Time Windows (CTOPTW), utilizing a Branch-and-Price scheme. To mitigate the computational burden inherent in exact methods, they effectively employed accelerated implicit enumeration and $ng$-route relaxations to solve the pricing subproblems. However, as the computational complexity of exact approaches often limits their applicability to smaller instances, subsequent research has increasingly favored heuristic strategies to balance solution quality with execution speed. Addressing this need for scalability, \cite{benSaid2019} proposed the Variable Space Search (VSS) heuristic. By alternating between giant tour and route search spaces within a hybrid framework of Greedy Randomized Adaptive Search Procedure (GRASP) and Evolutionary Local Search (ELS), this method offers a flexible mechanism to escape local optima that static search spaces might miss.

Building upon this heuristic foundation, \cite{hammami2024} recently demonstrated the superior efficacy of matheuristics by introducing the Hybrid Adaptive Large Neighborhood Search (HALNS). Unlike standard heuristics, HALNS integrates diverse removal and insertion operators with an exact Branch-and-Cut post-optimization step to solve a set packing formulation. This hybridization proved highly effective for large-scale benchmarks, identifying 10 new Best-Known Solutions (BKS) that prior methods failed to reach. While the aforementioned studies operate under deterministic assumptions, real-world logistics often involve significant uncertainty. Challenging the offline paradigm, \cite{shiri2024} investigated the Online CTOP (OCTOP), where parameters such as prizes, service times, and demands are revealed only upon arrival. They developed three online algorithms, deriving rigorous tight upper bounds on competitive ratios to ensure performance stability under uncertainty.

\subsubsection{Discussion}
Capacity constraints make OP variants particularly sensitive to the trade-off between rigor and scalability, and recent COP research reflects this tension. Early exact schemes such as the Branch-and-Price algorithm of \citet{park2017} clarify the structural underpinnings of CTOPTW, but also expose the size limits of purely exact optimization and thus motivate richer heuristic search spaces, exemplified by the variable-space design of \citet{benSaid2019}. Matheuristics like HALNS \citep{hammami2024} show how an exact set-packing component can be embedded within an adaptive large neighborhood search to recover much of the strength of exact methods on large instances, while the OCTOP model of \citet{shiri2024} pushes the field toward online, uncertainty-aware decision making. As a result, COP algorithms are increasingly evaluated not only by static solution quality but also by their responsiveness and robustness in dynamic operational settings.

\subsection{Clustered Orienteering Problem}\label{3.5}
The Clustered Orienteering Problem (CluOP) extends the classic OP by introducing clustering, where profits are tied to customer groups rather than individual locations.
This reflects logistics challenges (e.g., last-mile delivery), where all cluster members must be serviced to gain the profit. 
Since its introduction by \cite{angelelli2014}, the CluOP has inspired diverse computational advancements, ranging from exact methods to heuristics.
We review CluOP research, emphasizing methodological innovations for complex environments.

\subsubsection{Solution approaches}

\cite{yahiaoui2019} built upon the foundational work on the Clustered Orienteering Problem by \cite{angelelli2014}, extending the problem to incorporate multiple vehicles and formulating the Clustered Team Orienteering Problem.
The authors proposed an exact cutting-plane method and a hybrid heuristic combining Adaptive Large Neighborhood Search (ALNS) with a split procedure.

\cite{wu2024} proposed a Hybrid Evolutionary Algorithm (HEA) for the CluTOP.
This method integrates a backbone-based crossover operator for intensification and a destroy-and-repair mutation mechanism for diversity.
Additionally, a Reinforcement Learning (RL) guided Tabu Search refines exploration, focusing on high-potential candidates to reduce overhead.
\cite{he2024} proposed a Multi-Level Memetic Search, using a bilevel strategy combining cluster selection and routing.
\cite{montemanni2024} presented a Constraint Programming model to solve the CluTOP.

Recently, \cite{almeida2025} introduced the Clustered Orienteering Problem with Subgroups (CluOPS).
This variant hierarchically organizes nodes into subgroups within clusters, addressing the limitation of homogeneous clusters.
The model requires visiting all nodes within a subgroup to collect the reward, while restricting selection to at most one subgroup per cluster.
This framework theoretically unifies and generalizes both the classical CluOP and the Set OP (SOP).
The authors developed an Integer Linear Programming (ILP) formulation solved via Branch-and-Cut, and a Tabu Search metaheuristic for larger instances.

\subsubsection{Discussion}
Clustering fundamentally reshapes the OP landscape by shifting decisions from individual customers to groups, and recent CluOP studies increasingly exploit this structure rather than treating it as a mere side constraint. Algorithms such as cutting-plane–enhanced ALNS, hybrid evolutionary search with backbone operators and RL-guided tabu components, and multi-level memetic schemes show a clear move toward designs that intertwine cluster-aware representations with powerful local and global search mechanisms \citep{yahiaoui2019,wu2024,he2024}. The introduction of CluOPS \citep{almeida2025}, with its hierarchical subgroup structure and unifying view of CluOP and the Set OP, further underscores a shift from homogeneous clusters to richer intra-cluster heterogeneity. Taken together, these developments position clustered variants as a bridge between classical OPs and more general location–routing models, where the central challenge is to exploit multi-level group structure to balance profitability, coverage, and computational effort.

Table \ref{tab:12} in \ref{append2} summarizes the assumptions discussed in Sections \ref{3.1}–\ref{3.5}.
Table \ref{tab:13} details in \ref{append2} the solution methods covered in Sections \ref{3.1}–\ref{3.5}.

\subsection{Other Variants}\label{3.6}

\textbf{Synchronized Orienteering Problem} 

The Synchronized Orienteering Problem represents a class of cooperative OPs, requiring collaboration among multiple heterogeneous vehicles.

\cite{roozbeh2020} addressed the Cooperative Orienteering Problem with Time Windows (CoOPTW), an extension of the TOPTW.
Motivated by applications like emergency logistics and home healthcare, the authors addressed synchronization requirements at specific nodes, a feature absent in traditional models.
Using an Adaptive Large Neighborhood Search (ALNS) framework, they introduced merit-based heuristics, including tailored ruin-and-recreate strategies.
Addressing disaster management challenges, \cite{yahiaoui2023} focused on the Synchronized Team OP with Time Windows (STOPTW). This variant emphasizes temporal coordination and vehicle compatibility.
The authors employed a hybrid metaheuristic combining GRASP, Iterated Local Search (ILS), and set-covering optimization to handle large-scale instances.

\textbf{Arc Orienteering Problem} 

\cite{riera2017} addressed the Team Orienteering Arc Routing Problem (TOARP), where customers are represented by arcs in a directed graph.
The study was motivated by weak dual bounds in Linear Programming (LP) relaxations of knapsack constraints.
The authors proposed a column generation approach, transforming arc-based representation into vertex-based to simplify formulation and branching.

Recently, \cite{martin2025} extended the TOARP by introducing a variant allowing different origin and destination depots.
Motivated by UAV operations where drones deploy and recover at different locations, this model relaxes the single-depot assumption.
To address this, the authors proposed a Biased-Randomized Iterated Local Search (BR-ILS) algorithm.
This approach combines a biased-randomized savings heuristic with perturbation-based local search. It demonstrated high efficiency, outperforming Best-Known Solutions (BKS) on challenging benchmarks.

Table \ref{tab:14} summarizes the assumptions of articles in Section \ref{3.6}.
Table \ref{tab:15} details the solution methods for articles in Section \ref{3.6}.

\begin{table}[htbp]
	\centering
	\caption{Assumptions of other canonical extended models.}
	\scriptsize % 1. 使用小号字体 (如果还觉得大，可以改成 \tiny)
	\setlength{\tabcolsep}{1.5pt} % 2. 极大地减小列间距 (默认是 6pt)
	\renewcommand{\arraystretch}{0.85}
	\begin{tabular}{llcccc}
		\toprule
		Reference  &  Problem  & \multicolumn{1}{l}{Team} & \multicolumn{1}{l}{Time windows} & \multicolumn{1}{l}{Cooperative} & \multicolumn{1}{l}{Arc} \\
		\midrule
		\cite{roozbeh2020} & CoOP-TW & $\checkmark$ & \multicolumn{1}{c}{$\checkmark$} & \multicolumn{1}{c}{$\checkmark$} &  \\
		\cite{yahiaoui2023} & SynTOPTW & $\checkmark$ & \multicolumn{1}{c}{$\checkmark$} & \multicolumn{1}{c}{$\checkmark$} &  \\
		\cite{riera2017} & TOARP & $\checkmark$ &       &       & $\checkmark$ \\
		\cite{martin2025} & TOARP & $\checkmark$ & & & $\checkmark$ \\
		\bottomrule
	\end{tabular}%
	\label{tab:14}%
\end{table}%

\begin{table}[htbp]
	\centering
	\caption{Types of solution approaches for other canonical extended models.}
	\scriptsize % 1. 使用小号字体 (如果还觉得大，可以改成 \tiny)
	\setlength{\tabcolsep}{1pt} % 2. 极大地减小列间距 (默认是 6pt)
	\renewcommand{\arraystretch}{0.85}
	\begin{tabular}{llcccc}
		\toprule
		Reference  &  Algorithm & \multicolumn{1}{l}{Exact } & \multicolumn{1}{l}{Heuristic } & \multicolumn{1}{l}{Metaheuristic } & \multicolumn{1}{l}{Matheuristic }  \\
		\midrule
		\cite{roozbeh2020} & ALNS  &       &       & \multicolumn{1}{c}{$\checkmark$} & \\
		\cite{yahiaoui2023} & GRASP + ILS &       &       & \multicolumn{1}{c}{$\checkmark$} & \multicolumn{1}{c}{$\checkmark$}\\
		\cite{riera2017} & Column generation & $\checkmark$ &       &       & \\
		\cite{martin2025} & BR-ILS & & & $\checkmark$\\
		\bottomrule
	\end{tabular}%
	\label{tab:15}%
\end{table}

\subsection{Discussion}

The evolution of canonical OP extensions reveals a distinct shift from model definition to algorithmic sophistication. While early research primarily focused on formulating these variants, recent studies from 2017 to 2025 demonstrate a clear bifurcation in solution methodologies driven by problem characteristics.

First, for variants characterized by high uncertainty and dynamism (e.g., Stochastic and Time-Dependent OPs), there is a dominant trend toward integrating simulation with optimization. Approaches such as Simulation-Aided Multiple Plan Approaches (SMPA) \citep{bian2018} and Simheuristics \citep{panadero2022} have replaced simple re-optimization, reflecting a need for robust, anticipatory decision-making in logistics. Furthermore, the burgeoning integration of Machine Learning (e.g., Graph Convolutional Networks) is opening new frontiers, enabling real-time responsiveness in dynamic settings where traditional re-optimization may be computationally prohibitive.

Conversely, for variants defined by complex structural constraints (e.g., OPMV, CluOP), we observe a resurgence of exact methods. As heuristic feasibility becomes precarious under strict constraints like mandatory visits or synchronization, advanced Branch-and-Cut-and-Price algorithms \citep{peng2025, perez2025, almeida2025} and polyhedral studies have proven essential. These methods now scale to sizes previously reachable only by heuristics, challenging the historical ``heuristics-only'' assumption for complex OPs.

In summary, canonical OP extensions have matured into a rigorous testing ground for hybridizing OR and AI. They no longer merely represent ``added constraints'' but drive the development of specialized solvers that balance the trade-off between computational tractability and structural complexity.

\section{Emerging models and frontier extensions of the orienteering problem}\label{sec:4}
While the canonical extensions discussed in the previous section focus primarily on algorithmic refinements for established problems, the period from 2017 to 2025 has witnessed the emergence of novel OP variants that fundamentally reshape the problem's mathematical structure. These emerging models are driven by the need to capture increasingly complex real-world dynamics, such as the kinematic constraints of UAVs, the structural dependencies in last-mile delivery clusters, and the evolving nature of information in dynamic environments.

In this section, we apply the component-based taxonomy introduced in Section~\ref{sec:taxonomy} to systematically categorize these frontier extensions. Unlike the name-based organization of canonical variants, these emerging models are best understood through their structural alterations to the base OP. We structure the review as follows: Time-based emerging models (e.g., Multi-Period OP) that extend decision horizons; Node-based emerging models (e.g., Multi-Profit OP) that enrich objective functions; Path-based emerging models (e.g., Dubins OP) that introduce continuous motion constraints; Structure-based emerging models (e.g., Set OP) that modify topological dependencies; and Information-based emerging models (e.g., Probabilistic and Dynamic OP) that integrate learning and uncertainty management. This structured approach highlights how each component transforms the underlying optimization challenge, paving the way for the specialized solution strategies discussed therein.

Tables \ref{tab:18} and \ref{tab:20} in \ref{append3} summarize the modeling assumptions of the emerging new models discussed in Sections \ref{sec:4}. Tables \ref{tab:19} and \ref{tab:21} in \ref{append3} detail the corresponding solution methods.

\subsection{Time-based emerging models}
While canonical time-based extensions primarily focus on scheduling constraints within a single operational cycle (e.g., time windows), emerging models expand the temporal dimension to encompass multi-period planning horizons. In these variants, the optimization challenge shifts from purely spatial routing to a complex trade-off between immediate reward collection and long-term resource management across extended timelines. This structural evolution is critical for applications requiring sustained operations, such as infrastructure monitoring and security patrolling, where decisions in one period inherently constrain future possibilities.

\subsubsection{Multi-Period Orienteering Problem}

Multi-period variants optimize routes over extended horizons, balancing immediate rewards with long-term goals.

\textbf{Solution approaches.}
\cite{wang2024} applied the Multi-Period OP to post-earthquake building damage assessment. 
They formulated inspection routing as an MPOP to maximize information gain for a Gaussian Process Regression (GPR) model. A Lagrangian relaxation heuristic was employed to optimize routes across multiple days.
\cite{zhang2020} introduced a Multi-Period OP integrating stochastic wait times and dynamically evolving adoption likelihoods, modeling it as a Markov Decision Process (MDP).
\cite{vidigal2023} addressed a car patrolling problem for a security company. The study aimed to maximize the weighted sum of collected scores while minimizing patrol duration.

\subsubsection{Discussion}
Time-based emerging models thus extend classical temporal components by explicitly coupling decisions across multiple periods, and the existing multi-period OP studies consistently show that such long-horizon planning can outperform myopic period-by-period routing in settings ranging from post-disaster inspections to customer adoption and security patrols. 
At the same time, the literature in this area remains sparse and predominantly application-driven, with only a few structurally related formulations and no common benchmark set, suggesting that multi-period coupling is still treated as a niche modeling choice rather than a fully developed time-based component of the OP family. 
Future work could therefore focus on generic formulations and decomposition schemes that isolate inter-period interactions, on integrating multi-period decisions with stochastic or information-acquisition mechanisms, and on curated benchmark instances that clarify when the added complexity of multi-period models yields substantial value over repeated single-period orienteering.

\subsection{Node-based emerging models}
Unlike classical formulations that assign a single fixed score to each location, node-based emerging models enrich the problem structure by introducing complex reward attributes. These variants capture scenarios where the value of a visit is multifaceted, involving multiple objectives or time-dependent profitability. Consequently, the optimization challenge expands from simple routing to managing intricate trade-offs within the objective function itself.

\subsubsection{Multi-Profit Orienteering Problem}
This variant explores complex reward structures beyond simple static values, as each vertex is associated with multiple profit components that may depend on the visit time. The model thus explicitly captures heterogeneous stakeholder objectives (e.g., ticket sales versus publicity) within a single routing framework.

\textbf{Solution approaches.}
Kim and Kim (2020) introduced the Multi-Profit Orienteering Problem (MPOP), assigning multiple time-dependent profits to each vertex and proposing simulated annealing–based heuristics that adapt classical OP benchmark instances to the multi-profit setting. Building on this formulation, Kim and Kim (2022) developed the first exact algorithm based on hybrid Dynamic Programming with bounding, which solves many small and medium-size instances to optimality and provides strong reference bounds for larger cases.

\subsubsection{Discussion}
Node-based emerging models remain relatively scarce, indicating that complex node dynamics are still predominantly captured through rigid constraints rather than flexible, value-generating reward functions. Existing work shows that multi-attribute or time-dependent profits can better reflect service heterogeneity but also introduce significant modelling and algorithmic complexity. Future research could therefore focus on developing generic formulations and benchmark sets for multi-attribute reward structures and on exploring dynamic programming and decomposition-based approaches that remain computationally tractable.

\subsection{Path-based emerging models}
Path-based emerging models refine the path-based component of our taxonomy by replacing abstract arcs with motion primitives that reflect vehicle kinematics and geometric constraints. 
In contrast to classical formulations where travel costs are static attributes of arcs, these variants explicitly encode heading, turning-radius, energy, or neighborhood constraints in the edge definition. 
We group these developments into three main families: the Dubins Orienteering Problem, which imposes curvature constraints for nonholonomic vehicles; Drone and Robot Orienteering Problems, which incorporate kinematic limits, energy consumption, and geometric neighborhoods; and the Orienteering Problem with Drones, which coordinates ground vehicles and drones along coupled routes. 
Together, these models illustrate a shift from purely discrete graph routing toward hybrid formulations in which combinatorial route choices interact with (discretized) continuous motion planning.

\subsubsection{Dubins Orienteering Problem}\label{dop}
The Dubins Orienteering Problem (DOP) extends the classic OP by introducing curvature constraints for nonholonomic vehicles (e.g., UAVs and autonomous cars).
The DOP optimizes reward collection within a travel budget while adhering to turning radius constraints. 
This bridges theory and practice in applications like environmental monitoring.
\cite{penicka2017} formalized the DOP, integrating Dubins vehicle dynamics into path planning and target selection. 
Subsequent research advanced this with heuristics and metaheuristics, addressing computational challenges in balancing curvature constraints and path efficiency.
This section reviews recent DOP developments, emphasizing novel algorithms and real-world deployment.

\textbf{Solution approaches.}
\cite{penicka2017} introduced the DOP for curvature-constrained vehicles.
Addressing UAV data collection challenges, the study developed a Variable Neighborhood Search (VNS)
Building on the Correlated Orienteering Problem, \cite{tsiogkas2018} introduced the Dubins Correlated Orienteering Problem (DCOP). 
Motivated by resource-limited sensing missions in dynamic environments, the authors proposed a Genetic Algorithm heuristic to ensure near-optimal solutions while maintaining computational efficiency for real-time applications. 

\cite{faigl2020} tackled the Close Enough Dubins Orienteering Problem (CEDOP), allowing vehicles to collect rewards within a defined proximity of targets.
This formulation aims to reduce path lengths and complexity.
The study introduced a Growing Self-Organizing Array (GSOA) algorithm, integrating combinatorial and continuous optimization.
\cite{macharet2021} introduced the Minimal Exposure Dubins Orienteering Problem (MEDOP). 
This multi-objective formulation addresses reward maximization and exposure minimization.
Targeting surveillance and military operations, it employs an evolutionary algorithm to refine target selection, visitation order, and path curvature.

\subsubsection{Drone and Robot Orienteering Problem}

The term ``Drone and Robot Orienteering Problem'' originates from the integration of unmanned systems with the OP framework. 
Research typically adapts the model for kinematic constraints, geometric relaxations, or persistent monitoring.

\textbf{Solution approaches.}
\cite{wan2024} addressed a Multi-UAV scheduling problem for disaster data collection, where data value decays based on arrival and service time. 
Analyzing route selection and service time, they proposed an attention-based Deep Reinforcement Learning (DRL) framework for real-time decision-making.
Regarding geometric relaxations, \cite{qian2024} investigated the Close Enough Orienteering Problem (CEOP) with overlapped neighborhoods.
Allowing agents to collect rewards by entering a region, they used a Randomized Steiner Zone Discretization (RSZD) scheme to transform the problem into a Set OP (SOP).
Furthermore, \cite{asghar2024} applied OP concepts to multi-robot persistent monitoring. 
Although addressing a covering problem, the authors developed an orienteering-based heuristic to minimize the number of robots under latency and recharging constraints.

Regarding kinematic constraints, \cite{bayliss2020} introduced a learnheuristic for aerial drones, incorporating physical factors like air resistance and gravity.
\cite{sundar2022} addressed fixed-wing drone kinematics, employing a Branch-and-Price algorithm with heading angle discretization.

\subsubsection{Orienteering Problem with Drone} 
The Orienteering Problem with Drone addresses the coordinated routing of vehicles and drones.

\textbf{Solution approaches.}
\cite{morandi2024} introduced the Orienteering Problem with Multiple drones (OP-mD), featuring multiple drones collaborating with a truck.
They proposed a Mixed-Integer Linear Programming formulation and developed a tailored Branch-and-Cut algorithm.
Experimental results demonstrated the algorithm's effectiveness, solving instances with up to 50 nodes. 
The authors adapted their framework to address related problems, comparing computational performance against previous studies.

\subsubsection{Discussion}
Path-based emerging models have proliferated in recent years, indicating a shift from abstract graph routing toward formulations that explicitly encode kinematic, energy, and neighborhood constraints in UAV and mobile-robot applications. Existing studies combine advanced decomposition methods with tailored heuristics and, more recently, learning-augmented approaches, but they still rely on problem-specific discretizations and lack shared benchmark libraries across Dubins, close-enough, and truck–drone variants. Future work may therefore seek more unified models and instance sets that bridge discrete route choice and continuous trajectory planning, clarifying when the additional geometric fidelity of path-based formulations justifies their computational cost.

\subsection{Structure-based emerging models}
Structure-based emerging models modify the combinatorial structure of the vertex set rather than adding local constraints or temporal resources. Instead of deciding on individual customers, these variants operate on sets, clusters, or interacting routes, capturing group-level coverage requirements and coordination effects that arise in logistics, distribution, and multi-agent planning. This shift aligns OP research with generalized and clustered routing models and naturally motivates set-partitioning formulations and branch-and-price or matheuristic frameworks. In what follows, we focus on the Set Orienteering Problem and its recent generalizations as a representative class of structure-based emerging models.

\subsubsection{Set Orienteering Problem}\label{sop}
Recently, the Set Orienteering Problem (SOP) has attracted significant research attention.
Introduced by \cite{archetti2018}, the SOP generalizes the OP by grouping customers into clusters, each assigned a profit.
Profit is collected if the vehicle visits at least one customer within the cluster.
The objective is to design a single-vehicle route that maximizes total profit without exceeding a predefined duration threshold.
The SOP applies to logistics and distribution systems.
For instance, visiting a representative customer in a regional cluster may satisfy the entire demand. This approach optimizes costs and improves efficiency.
Table \ref{tab:SOP} lists the benchmark instances for S(T)OP(TW).

\begin{table}[htbp]
	\centering
	\scriptsize % 1. 使用小号字体 (如果还觉得大，可以改成 \tiny)
	\setlength{\tabcolsep}{1.5pt} % 2. 极大地减小列间距 (默认是 6pt)
	\renewcommand{\arraystretch}{0.85}
	\caption{Benchmark instances of S(T)OP(TW).}
	\begin{tabular}{crr}
		\toprule
		\multicolumn{1}{l}{Reference} & \multicolumn{1}{l}{Number of Instances} & \multicolumn{1}{l}{Number of Nodes ($|$N$|$)} \\
		\midrule
		\multirow{5}[1]{*}{\cite{archetti2018}} & 6     & 52 \\
		& 6     & 51 \\
		& 4     & 70 \\
		& 4     & 76 \\
		& 24    & 200 to 1084 \\
		\cite{dontas2023} & 19$\times$16 & 1000 to 3162 \\
		\bottomrule
	\end{tabular}%
	\label{tab:SOP}%
\end{table}%

\textbf{Solution approaches.}
As a foundational study on the Set Orienteering Problem, \cite{archetti2018} introduced the SOP as a generalization of the OP and proposed the MASOP algorithm.
This matheuristic combines a greedy construction phase with Tabu Search for refinement.
The study established a foundation for SOP solutions and demonstrated their logistics applicability.

\cite{penicka2019} applied Variable Neighborhood Search (VNS) to the SOP, integrating Integer Linear Programming (ILP) to solve small-scale instances.
This study emphasized the versatility of VNS across SOP variants. 
\cite{carrabs2021} introduced the Biased Random-Key Genetic Algorithm (BRKGA).
Using pre-processing and local search, BRKGA outperformed VNS and MA-SOP in both computation and quality on standard benchmarks.
\cite{dontas2023} developed an Adaptive Memory Matheuristic (AMM) that integrates local search with mathematical programming. 
AMM significantly outperformed MASOP, consistently surpassing \cite{archetti2018} benchmarks, particularly on large datasets.
This highlights the role of adaptive memory in generating diverse, high-quality solutions.
\cite{lu2024} presented a Hybrid Evolutionary Algorithm (HEA) that combines cluster-based crossover, randomized mutation, and local refinement procedures.

\cite{archetti2024} introduced a novel formulation and a Branch-and-Cut algorithm.
By reducing variable dimensionality and adding valid inequalities, this method outperformed previous exact approaches on instances with up to 100 customers.
\cite{dutta2020} developed a multi-objective SOP model balancing customer satisfaction and profit maximization, adding flexibility through third-party logistics. 
The authors employed NSGA-II and SPEA2 for optimization, achieving improved results on synthesized datasets. 
\cite{lin2024} addressed the Set Orienteering Problem with Mandatory Visits (SOPMV), combining Simulated Annealing (SA) with a two-stage MILP model.
This method handles mandatory constraints, extending the SOP to scenarios requiring specific visits.
\cite{yu2024} developed a hybrid Simulated Annealing and Reinforcement Learning (SARL) approach for the Set Team OP with Time Windows (STOPTW).
This approach surpassed traditional algorithms in quality and efficiency, demonstrating adaptability to new variants.

Most recently, \cite{nguyen2025} generalized the SOP to the multi-vehicle Set Team Orienteering Problem (STOP).
They formulated a MILP model and developed a Branch-and-Price (B\&P) algorithm incorporating bucket-based dynamic programming and strong branching.
To handle large-scale instances, they proposed a Large Neighborhood Search (LNS) metaheuristic with novel proximity-based insertion operators.
Notably, these contributions extend to standard SOP benchmarks. 
The B\&P algorithm solved $61.9\%$ of instances to optimality, outperforming previous exact methods. 
Additionally, the LNS discovered 40 new Best-Known Solutions (BKS).

\subsubsection{Discussion}
Structure-based emerging models for the OP have grown in recent years, but the literature is still concentrated on set- and cluster-based formulations such as the SOP and its generalizations. Existing studies show that group-level routing decisions can be handled effectively by modern exact, matheuristic, and learning-augmented approaches, and that there is still a lack of systematic benchmarks across structural variants. Future work may develop scalable algorithms for team-based structural variants and explore richer topological dependencies.

\subsection{Information-based emerging models}
Information-based emerging models extend the classical Stochastic Orienteering Problem (StchOP) by addressing uncertainty through lenses beyond standard risk-neutral expectation. While StchOP typically optimizes expected values under known probability distributions for travel times or rewards, the variants in this section grapple with the quality, availability, and evolution of information. Specifically, the Probabilistic OP shifts focus to the stochastic availability of customers, whereas Uncertain and Robust OPs tackle parameter ambiguity and adversarial worst-case scenarios without relying on precise distributional assumptions. Additionally, Dynamic OPs take into account the issue of information disclosure, requiring the development of adaptive strategies based on the information revealed during the execution process.

\subsubsection{Probabilistic Orienteering Problem}\label{pop}

The Probabilistic Orienteering Problem (POP) is a stochastic extension that differs from the Stochastic OP (StchOP), which typically assumes randomness in rewards or travel times.
Instead, the POP focuses on customer availability uncertainty, capturing real-world decision-making complexities.
Optimizing the trade-off between expected rewards and costs, the POP addresses challenges in dynamic logistics and adaptive planning.
\cite{angelelli2017} established the POP foundation with a mathematical formulation and solution strategies, including Branch-and-Cut (B\&C) algorithms.
Subsequent research employs advanced metaheuristics and simulation to address complex scenarios.

\textbf{Solution approaches.}
\cite{angelelli2017} conceptualized the POP as a two-stage decision-making framework.
First, a subset of nodes is selected, and a preliminary path is planned. 
Then, the path is adjusted dynamically once node availability is revealed.
Using a stochastic MILP formulation, the authors developed a Branch-and-Cut algorithm and matheuristics to improve efficiency.
The goal is to maximize expected profits (collected rewards minus costs).
\cite{chou2020} introduced a Tabu Search algorithm with Monte Carlo sampling for solution evaluation.
They optimized prize collection under uncertainty by embedding stochastic approximation into the metaheuristic.
This approach achieved both precision and speed, escaping local optima.

\cite{angelelli2021} introduced a dynamic and probabilistic variant (DPOP), integrating real-time decision-making.
Framed as a Markov Decision Process (MDP), it enables optimal policies maximizing expected profits.
To address computational challenges, the authors designed heuristics, including Sample Average Approximation (SAA).
\cite{herrera2022} investigated the Team Orienteering Problem with Probabilistic Delays (TOP-PD), addressing travel time uncertainties.
Modeling delays as probabilistic distributions, the study incorporated a reliability constraint to enforce delay thresholds.
The authors proposed a simheuristic algorithm using reliability analysis to evaluate solutions via survival functions.

Recently, \cite{montemanni2025} proposed an exact algorithm iteratively solving deterministic approximations.
They introduced two deterministic models to bound expected travel times, using Constraint Programming (CP-SAT) to find optimal solutions.
Experiments showed this approach outperforms previous methods in solution quality and speed on standard benchmarks.

\subsubsection{Dynamic Orienteering Problem}

The Dynamic Orienteering Problem (DyOP) extends the classic OP by incorporating dynamic environmental changes, such as varying customer availability.
This model optimizes decision-making amidst evolving information, requiring adaptive strategies to maximize rewards.

\textbf{Solution approaches.}
Recent research has formalized the definitions and expanded the decision scope of the DyOP. 
\cite{kirac2025} formally defined the Dynamic Team Orienteering Problem (DyTOP), where customer locations are revealed over time (deterministic dynamics).
They adapted a Multiple Plan Approach (MPA) to maintain a pool of routing plans, using consensus functions to select actions.
Focusing on service logistics, \cite{paradiso2025} addressed a dynamic next-day service problem involving hard time windows for dynamic requests.
They proposed a stochastic lookahead algorithm using a linear relaxation of a set packing reformulation. This evaluates the future flexibility of current assignments.

\cite{wang2023} addressed a dynamic and stochastic OP variant in autonomous transportation systems.
The study investigates real-time routing challenges, where stochastic times require adaptive strategies.
The authors introduced three multi-stage heuristics, including Insertion-and-Removal and Hybrid Simulated Annealing-Tabu Search Re-Optimization (HSATSRO).
\cite{saeedvand2020} introduced a Hybrid Multi-Objective Algorithm for rescue operations, combining Evolutionary Algorithms with Q-learning.
The approach uses Reinforcement Learning to dynamically adapt to environmental changes, significantly outperforming traditional methods. 
\cite{uguina2024} introduced a heterogeneous, dynamic Team OP with probabilistic rewards influenced by environmental factors.
They proposed a learnheuristic framework combining Thompson sampling with metaheuristics.
\cite{li2024} investigated sequential decision-making in same-day delivery, modeling it as the OP with Stochastic and Dynamic Release Dates and using Reinforcement Learning.

Furthermore, some studies have introduced highly interdependent dynamic constraints. 
\cite{granda2025} proposed the Team Orienteering Problem with Variable Time Windows (TOPVTW), motivated by wildfire suppression. 
Unlike problems driven by external information, this model features endogenous dynamics: visiting a node ``blocks'' fire spread, dynamically extending other time windows.

\subsubsection{Uncertain Orienteering Problem}

The Uncertain Orienteering Problem integrates uncertainty theory into the classic OP, addressing imprecise or unpredictable factors.

\textbf{Solution approaches.}
Moving beyond models with known distributions, \cite{qian2025} proposed the Uncertain and Dynamic Orienteering Problem (UDyOP).
Their model assumes edge costs follow distributions with unknown and time-variant parameters. 
To solve this, they developed an Adaptive Approach for Probabilistic Paths (ADAPT), using Bayesian inference to update beliefs and estimate costs online.

Earlier works focused on different aspects of uncertainty. 
\cite{wang2019a} integrated uncertainty theory into the OP, introducing a minimum-risk model based on uncertain measures rather than expected values.
The study formulated an OP with uncertain profits and proposed a One-Dimensional Ratio Search algorithm.
\cite{jin2019} addressed routing phlebotomists, where rewards and service times are stochastic.
The authors developed an \textit{a priori} policy using Monte Carlo simulations to estimate solution values.
\cite{wang2019b} extended the framework to include team dynamics and time windows, formulating the Uncertain Team OP with Time Windows (UTOPTW).

\subsubsection{Robust Orienteering Problem}

The Robust Orienteering Problem (ROP) incorporates robustness against uncertainties (e.g., travel time variations) to ensure reliable solutions under worst-case scenarios.

\textbf{Solution approaches.}
\cite{yu2022} addressed a robust Team OP with Decreasing Profits, where profits decay and service times are uncertain.
They proposed a robust optimization framework incorporating Branch-and-Price and Tabu Search.
\cite{zhang2023} proposed a robust drone arc routing framework targeting humanitarian logistics. 
The model addresses road condition uncertainty, using graph transformation to convert the multigraph problem into a node-based formulation.
\cite{shi2023} introduced the Robust Multiple-Path Orienteering Problem (RMOP) to address the vulnerability of robot teams in adversarial environments. 
For the offline problem, they used approximation algorithms with modular/submodular rewards. For the online problem, they introduced an MCTS-based strategy.

\subsubsection{Discussion}
Information-based emerging models show how uncertainty, information revelation, and adversarial conditions reshape OP decisions beyond risk-neutral expectations. Recent work on probabilistic, dynamic, uncertain, and robust variants shifts the focus from fixed routes to policies and learning-enhanced heuristics that react to incoming information and manage risk. These developments point to future research that more tightly links data-driven learning, robustness, and multi-stage decision frameworks.

\subsection{Discussion}
Recent OP extensions, including the Set, Dubins, and Probabilistic Orienteering Problems, address complex real-world challenges.
These variants incorporate uncertainty, dynamic conditions, and specialized constraints. 
This makes them suitable for logistics, robotics, and disaster response, where traditional assumptions are insufficient.
Alongside Dynamic and Drone OPs, they enhance relevance in stochastic environments and settings like autonomous transportation.

Solution methods often leverage hybrid and metaheuristic techniques (e.g., Tabu Search, Genetic Algorithms) to tackle these NP-hard problems efficiently.
Advanced approaches (e.g., Reinforcement Learning, Monte Carlo simulation) address uncertainty, while exact methods (e.g., Branch-and-Cut) suit smaller instances.
This reflects a shift toward robust optimization frameworks, driven by innovative modeling and sophisticated strategies.

\section{Conclusion and future research directions}\label{sec5}

Since \cite{tsiligirides1984} introduced the OP, it has evolved from a basic path optimization challenge into a flexible framework addressing complex real-world scenarios.
This review tracks the expansion of OP models—from single-person variants to the Team OP, OP with Time Windows, and modern Dynamic, Stochastic, and Set OPs—reflecting a growing demand for adaptable tools under uncertainty.
Applications have broadened to logistics, tourism, urban management, UAV navigation, and disaster response.
As models advance, so do solution methods.
Traditional exact methods have improved significantly; notably, \cite{kobeaga2024} redesigned a Branch-and-Cut algorithm to tackle instances with up to 7,397 nodes.
Meanwhile, heuristics (e.g., ALNS) and Machine Learning (ML) techniques, particularly Reinforcement Learning (RL) and neural networks, have enhanced scalability and real-time decision-making for dynamic variants.
These developments demonstrate a synergy between exact and approximate paradigms, balancing efficiency with solution quality.
Current research indicates significant progress, handling previously intractable instances and providing robust solutions that translate into tangible benefits. 
Examples include enhanced logistics efficiency and stronger emergency response capabilities.
Driven by innovative models and algorithms, the OP has matured into a cornerstone for solving resource-constrained optimization challenges.

\subsection{Future research directions}\label{sec6}

Despite remarkable achievements, problem scale and complexity necessitate further innovation. 
To guide future work, we highlight four priority avenues in Table \ref{tab:future_directions_compact}.

\begin{table}[htbp]
	\centering
	\scriptsize % 使用小号字体以节省空间
	\caption{Summary of future research directions.}
	\label{tab:future_directions_compact}
	\setlength{\tabcolsep}{4pt} % 减小列间距
	\renewcommand{\arraystretch}{1.1} % 稍微调整行间距，保持紧凑但可读
	\begin{tabular}{@{} p{0.28\linewidth} p{0.68\linewidth} @{}} % @{} 去除表格左右边缘的留白以最大化利用宽度
		\toprule
		\textbf{Theme} & \textbf{Key Research Avenues} \\
		\midrule
		\textbf{AI Integration} & Moving beyond basic RL to Generative AI/LLMs for high-level heuristics and Explainable AI (XAI) for trusted automation. \\
		\addlinespace[4pt] % 增加一点段落间距，比画横线更省空间且美观
		\textbf{Sustainability} & Developing Green OPs that balance profit with carbon footprints using realistic energy models. \\
		\addlinespace[4pt]
		\textbf{Dynamics \& Robustness} & Leveraging Distributionally Robust Optimization (DRO) and Digital Twins for ambiguity and real-time adaptation. \\
		\addlinespace[4pt]
		\textbf{Scalability} & Enhancing exact solvers via ML-aided Branching and advancing decomposition for multi-stage stochastic models. \\
		\bottomrule
	\end{tabular}
\end{table}

\textbf{Scalability and AI Integration.} 
Future research should enhance exact methods through parallel computing and ML-aided Branch-and-Bound, where models learn branching policies to accelerate convergence.
Simultaneously, metaheuristics can benefit from Generative AI and Large Language Models (LLMs) acting as high-level constructive heuristics.
While initial studies leveraged basic RL, recent breakthroughs like TOP-Former \citep{fuertes2025} demonstrate the power of Transformer-based architectures for multi-agent coordination.
Future work should leverage such architectures to predict optimal paths or learn from historical data \citep{li2024}.
Furthermore, as AI decision-making becomes more complex, developing Explainable AI (XAI) models is crucial to build trust in automated routing systems.

\textbf{Robustness and Uncertainty.} 
Real-world volatility necessitates resilience. 
Beyond standard stochastic methods, future research should explore Distributionally Robust Optimization (DRO) to handle ambiguity in probability distributions \cite{carpin2024, yu2022}.
Additionally, the integration of Digital Twins using IoT data streams offers a promising path for instant dynamic re-routing in hyper-dynamic environments.
Researching multi-objective variants and identifying Pareto-optimal solutions \cite{dasdemir2022} will further provide flexible decision-making tools for conflicting objectives.

\textbf{Sustainability and Green Logistics.} 
Research must transcend traditional profit maximization to address the ``Green Orienteering Problem.''
Current works like the Energy-Constrained OP (ECOP) \cite{karabas2025} utilize the Comprehensive Modal Emission Model (CMEM) to reduce GHG emissions.
Future frontiers include modeling non-linear battery degradation and integrating charging station planning into OP variants, addressing the unique constraints of electric fleets and drones.

\textbf{Domain-specific Applications.} 
Tailoring OP models to specialized domains strengthens practical relevance.
Incorporating kinematic constraints for UAVs \cite{bayliss2020, zhang2023} remains a key area.
Moreover, the OP framework shows potential in non-traditional fields, such as bioinformatics for DNA string packing \cite{andreani2024} and public health for disinfection robot scheduling \cite{liu2024}.
Deepening theoretical understanding—such as complexity analysis and polyhedral studies—remains essential to guide these algorithmic developments.

In summary, the future of OP lies in interdisciplinary efforts. Balancing computational cost with solution effectiveness, this field will continue to provide powerful tools for complex real-world challenges.

\clearpage
\appendix
	
	\section{New best known solution}\label{append1}
	
	\begin{sidewaystable}
		\centering
		\caption{Generation 1, n $>$ 400 and Generation 2, n $>$ 400.}
		\scriptsize % 1. 使用小号字体 (如果还觉得大，可以改成 \tiny)
		\setlength{\tabcolsep}{4pt} % 2. 极大地减小列间距 (默认是 6pt)
		\renewcommand{\arraystretch}{1}
		\begin{tabular}{rrrrrrrllllll}
			\cmidrule{1-6}\cmidrule{8-13}    \multicolumn{1}{l}{generation and size} & \multicolumn{1}{l}{instance} & \multicolumn{1}{l}{EA4OP} & \multicolumn{1}{l}{ALNS} & \multicolumn{1}{l}{RB\&C} & \multicolumn{1}{l}{BKs} &       & generation and size & instance & EA4OP & ALNS  & RB\&C & BKs \\
			\cmidrule{1-6}\cmidrule{8-13}    \multicolumn{1}{c}{\multirow{27}[2]{*}{gen1-large}} & \multicolumn{1}{l}{ali535} &       &       & \multicolumn{1}{l}{$\checkmark$} & 425   &       & \multicolumn{1}{c}{\multirow{29}[4]{*}{gen2-large}} & p654  &       & $\checkmark$ &       & \multicolumn{1}{r}{17900} \\
			& \multicolumn{1}{l}{pa561} &       &       & \multicolumn{1}{l}{$\checkmark$} & 357*  &       &       & gr666 &       &       & $\checkmark$ & \multicolumn{1}{r}{26514} \\
			& \multicolumn{1}{l}{dsj1000} &       &       & \multicolumn{1}{l}{$\checkmark$} & 656*  &       &       & rat783 &       &       & $\checkmark$ & \multicolumn{1}{r}{25474*} \\
			& \multicolumn{1}{l}{pr1002} &       &       & \multicolumn{1}{l}{$\checkmark$} & 606*  &       &       & dsj1000 &       &       & $\checkmark$ & \multicolumn{1}{r}{35835} \\
			& \multicolumn{1}{l}{u1060} &       &       & \multicolumn{1}{l}{$\checkmark$} & 660*  &       &       & pr1002 &       &       & $\checkmark$ & \multicolumn{1}{r}{33030} \\
			& \multicolumn{1}{l}{pcb1173} &       &       & \multicolumn{1}{l}{$\checkmark$} & 675*  &       &       & u1060 &       &       & $\checkmark$ & \multicolumn{1}{r}{36151} \\
			& \multicolumn{1}{l}{d1291} &       &       & \multicolumn{1}{l}{$\checkmark$} & 715*  &       &       & vm1084 &       &       & $\checkmark$ & \multicolumn{1}{r}{40777} \\
			& \multicolumn{1}{l}{rl1304} &       &       & \multicolumn{1}{l}{$\checkmark$} & 802*  &       &       & pcb1173 &       &       & $\checkmark$ & \multicolumn{1}{r}{37035} \\
			& \multicolumn{1}{l}{rl1323} &       &       & \multicolumn{1}{l}{$\checkmark$} & 814*  &       &       & d1291 &       &       & $\checkmark$ & \multicolumn{1}{r}{37778} \\
			& \multicolumn{1}{l}{nrw1379} &       &       & \multicolumn{1}{l}{$\checkmark$} & 815   &       &       & rl1304 &       &       & $\checkmark$ & \multicolumn{1}{r}{42275} \\
			& \multicolumn{1}{l}{fl1400} &       & \multicolumn{1}{l}{$\checkmark$} &       & 1048  &       &       & rl1323 &       &       & $\checkmark$ & \multicolumn{1}{r}{43377} \\
			& \multicolumn{1}{l}{u1432} &       &       & \multicolumn{1}{l}{$\checkmark$} & 754   &       &       & nrw1379 &       &       & $\checkmark$ & \multicolumn{1}{r}{46676} \\
			& \multicolumn{1}{l}{fl1577} &       &       & \multicolumn{1}{l}{$\checkmark$} & 897   &       &       & ﬂ1400 &       & $\checkmark$ &       & \multicolumn{1}{r}{56692} \\
			& \multicolumn{1}{l}{d1655} &       &       & \multicolumn{1}{l}{$\checkmark$} & 922   &       &       & u1432 &       &       & $\checkmark$ & \multicolumn{1}{r}{46946} \\
			& \multicolumn{1}{l}{vm1748} &       &       & \multicolumn{1}{l}{$\checkmark$} & 1276  &       &       & ﬂ1577 & $\checkmark$ &       &       & \multicolumn{1}{r}{45505} \\
			& \multicolumn{1}{l}{u1817} &       &       & \multicolumn{1}{l}{$\checkmark$} & 983*  &       &       & d1655 & $\checkmark$ &       &       & \multicolumn{1}{r}{47211} \\
			& \multicolumn{1}{l}{rl1889} &       &       & \multicolumn{1}{l}{$\checkmark$} & 1226* &       &       & vm1748 &       &       & $\checkmark$ & \multicolumn{1}{r}{68042} \\
			& \multicolumn{1}{l}{d2103} &       &       & \multicolumn{1}{l}{$\checkmark$} & 1200* &       &       & u1817 &       &       & $\checkmark$ & \multicolumn{1}{r}{54245} \\
			& \multicolumn{1}{l}{u2152} &       &       & \multicolumn{1}{l}{$\checkmark$} & 1151* &       &       & rl1889 &       &       & $\checkmark$ & \multicolumn{1}{r}{63308} \\
			& \multicolumn{1}{l}{u2319} & \multicolumn{1}{l}{$\checkmark$} &       &       & 1167  &       &       & d2103 &       &       & $\checkmark$ & \multicolumn{1}{r}{63426*} \\
			& \multicolumn{1}{l}{pr2392} &       &       & \multicolumn{1}{l}{$\checkmark$} & 1316  &       &       & u2152 &       &       & $\checkmark$ & \multicolumn{1}{r}{64649} \\
			& \multicolumn{1}{l}{pcb3038} &       &       & \multicolumn{1}{l}{$\checkmark$} & 1727  &       &       & u2319 &       &       & $\checkmark$ & \multicolumn{1}{r}{80914} \\
			& \multicolumn{1}{l}{fl3795} &       &       & \multicolumn{1}{l}{$\checkmark$} & 1965  &       &       & pr2392 &       &       & $\checkmark$ & \multicolumn{1}{r}{72843} \\
			& \multicolumn{1}{l}{fnl4461} &       &       & \multicolumn{1}{l}{$\checkmark$} & 2541  &       &       & pcb3038 &       &       & $\checkmark$ & \multicolumn{1}{r}{97902} \\
			& \multicolumn{1}{l}{rl5915} &       &       & \multicolumn{1}{l}{$\checkmark$} & 3593  &       &       & ﬂ3795 & $\checkmark$ &       &       & \multicolumn{1}{r}{103397} \\
			& \multicolumn{1}{l}{rl5934} &       &       & \multicolumn{1}{l}{$\checkmark$} & 3632  &       &       & fnl4461 &       &       & $\checkmark$ & \multicolumn{1}{r}{147109} \\
			& \multicolumn{1}{l}{pla7397} &       &       & \multicolumn{1}{l}{$\checkmark$} & 5289  &       &       & rl5915 &       &       & $\checkmark$ & \multicolumn{1}{r}{184424} \\
			\cmidrule{2-6}    \multicolumn{1}{c}{\multirow{3}[4]{*}{gen2-large}} & \multicolumn{1}{l}{ﬂ417} &       &       & \multicolumn{1}{l}{$\checkmark$} & 11933 &       &       & rl5934 &       &       & $\checkmark$ & \multicolumn{1}{r}{187034} \\
			& \multicolumn{1}{l}{d493} &       &       & \multicolumn{1}{l}{$\checkmark$} & 16995 &       &       & pla7397 &       &       & $\checkmark$ & \multicolumn{1}{r}{281977} \\
			\cmidrule{8-13}          & \multicolumn{1}{l}{att532} &       &       & \multicolumn{1}{l}{$\checkmark$} & 19635 &       & \multicolumn{6}{l}{*Optimal solutions obtained by exact algorithms} \\
			\cmidrule{1-6}          &       &       &       &       &       &       & \multicolumn{6}{l}{**Verified as optimal via exact methods} \\
		\end{tabular}%
		\label{tab:op1}%
	\end{sidewaystable}
	
	\begin{sidewaystable}
		\centering
		\caption{Generation 3, n $>$ 400 and Generation 4, n $>$ 400.}
		\scriptsize % 1. 使用小号字体 (如果还觉得大，可以改成 \tiny)
		\setlength{\tabcolsep}{4pt} % 2. 极大地减小列间距 (默认是 6pt)
		\renewcommand{\arraystretch}{1}
		\begin{tabular}{rrrrrrrllllll}
			\cmidrule{1-6}\cmidrule{8-13}    \multicolumn{1}{l}{generation and size} & \multicolumn{1}{l}{instance} & \multicolumn{1}{l}{EA4OP} & \multicolumn{1}{l}{ALNS} & \multicolumn{1}{l}{RB\&C} & \multicolumn{1}{l}{BKs} &       & generation and size & instance & EA4OP & ALNS  & RB\&C & BKs \\
			\cmidrule{1-6}\cmidrule{8-13}    \multicolumn{1}{c}{\multirow{25}[4]{*}{gen3-large}} & \multicolumn{1}{l}{pr439} &       &       & \multicolumn{1}{l}{$\checkmark$} & 15176 &       & \multicolumn{1}{c}{\multirow{4}[2]{*}{gen3-large}} & fnl4461 &       &       & $\checkmark$ & \multicolumn{1}{r}{146995} \\
			& \multicolumn{1}{l}{pcb442} &       &       & \multicolumn{1}{l}{$\checkmark$} & 14819** &       &       & rl5915 &       &       & $\checkmark$ & \multicolumn{1}{r}{203695} \\
			& \multicolumn{1}{l}{ali535} &       &       & \multicolumn{1}{l}{$\checkmark$} & 9414  &       &       & rl5934 &       &       & $\checkmark$ & \multicolumn{1}{r}{212021} \\
			& \multicolumn{1}{l}{gr666} &       &       & \multicolumn{1}{l}{$\checkmark$} & 17023 &       &       & pla7397 &       &       & $\checkmark$ & \multicolumn{1}{r}{322285} \\
			\cmidrule{9-13}          & \multicolumn{1}{l}{dsj1000} &       &       & \multicolumn{1}{l}{$\checkmark$} & 31434 &       & \multicolumn{1}{c}{\multirow{21}[2]{*}{gen4-large}} & pcb442 &       & $\checkmark$ &       & \multicolumn{1}{r}{5869} \\
			& \multicolumn{1}{l}{pr1002} &       &       & \multicolumn{1}{l}{$\checkmark$} & 39526* &       &       & pa561 & $\checkmark$ &       &       & \multicolumn{1}{r}{27719} \\
			& \multicolumn{1}{l}{u1060} &       &       & \multicolumn{1}{l}{$\checkmark$} & 37492 &       &       & u1060 & $\checkmark$ &       &       & \multicolumn{1}{r}{51775} \\
			& \multicolumn{1}{l}{vm1084} &       &       & \multicolumn{1}{l}{$\checkmark$} & 37669* &       &       & pcb1173 & $\checkmark$ &       &       & \multicolumn{1}{r}{56010} \\
			& \multicolumn{1}{l}{pcb1173} &       &       & \multicolumn{1}{l}{$\checkmark$} & 41257* &       &       & rl1323 & $\checkmark$ &       &       & \multicolumn{1}{r}{65664} \\
			& \multicolumn{1}{l}{d1291} &       &       & \multicolumn{1}{l}{$\checkmark$} & 41509 &       &       & nrw1379 & $\checkmark$ &       &       & \multicolumn{1}{r}{69214} \\
			& \multicolumn{1}{l}{rl1304} &       &       & \multicolumn{1}{l}{$\checkmark$} & 41881 &       &       & ﬂ1400 &       & $\checkmark$ &       & \multicolumn{1}{r}{70511} \\
			& \multicolumn{1}{l}{rl1323} &       &       & \multicolumn{1}{l}{$\checkmark$} & 47213 &       &       & ﬂ1577 & $\checkmark$ &       &       & \multicolumn{1}{r}{33754} \\
			& \multicolumn{1}{l}{nrw1379} &       &       & \multicolumn{1}{l}{$\checkmark$} & 42920 &       &       & d1655 &       & $\checkmark$ &       & \multicolumn{1}{r}{33231} \\
			& \multicolumn{1}{l}{ﬂ1400} &       & \multicolumn{1}{l}{$\checkmark$} &       & 57470 &       &       & vm1748 & $\checkmark$ &       &       & \multicolumn{1}{r}{82126} \\
			& \multicolumn{1}{l}{u1432} &       &       & \multicolumn{1}{l}{$\checkmark$} & 47778 &       &       & u1817 &       & $\checkmark$ &       & \multicolumn{1}{r}{37457} \\
			& \multicolumn{1}{l}{ﬂ1577} & \multicolumn{1}{l}{$\checkmark$} &       &       & 45692 &       &       & rl1889 &       & $\checkmark$ &       & \multicolumn{1}{r}{83875} \\
			& \multicolumn{1}{l}{d1655} &       &       & \multicolumn{1}{l}{$\checkmark$} & 62048 &       &       & d2103 &       & $\checkmark$ &       & \multicolumn{1}{r}{37124} \\
			& \multicolumn{1}{l}{vm1748} &       &       & \multicolumn{1}{l}{$\checkmark$} & 71885 &       &       & u2152 &       & $\checkmark$ &       & \multicolumn{1}{r}{55397} \\
			& \multicolumn{1}{l}{u1817} & \multicolumn{1}{l}{$\checkmark$} &       &       & 63639 &       &       & pr2392 &       & $\checkmark$ &       & \multicolumn{1}{r}{50944} \\
			& \multicolumn{1}{l}{rl1889} &       &       & \multicolumn{1}{l}{$\checkmark$} & 70065 &       &       & pcb3038 & $\checkmark$ &       &       & \multicolumn{1}{r}{101173} \\
			& \multicolumn{1}{l}{d2103} &       &       & \multicolumn{1}{l}{$\checkmark$} & 82787 &       &       & ﬂ3795 & $\checkmark$ &       &       & \multicolumn{1}{r}{80069} \\
			& \multicolumn{1}{l}{u2152} &       &       & \multicolumn{1}{l}{$\checkmark$} & 74007 &       &       & fnl4461 & $\checkmark$ &       &       & \multicolumn{1}{r}{85088} \\
			& \multicolumn{1}{l}{pr2392} &       &       & \multicolumn{1}{l}{$\checkmark$} & 85409 &       &       & rl5915 & $\checkmark$ &       &       & \multicolumn{1}{r}{279277} \\
			& \multicolumn{1}{l}{pcb3038} &       &       & \multicolumn{1}{l}{$\checkmark$} & 106928 &       &       & rl5934 & $\checkmark$ &       &       & \multicolumn{1}{r}{137838} \\
			& \multicolumn{1}{l}{ﬂ3795} & \multicolumn{1}{l}{$\checkmark$} &       &       & 97707 &       &       & pla7397 & $\checkmark$ &       &       & \multicolumn{1}{r}{142399} \\
			\cmidrule{1-6}\cmidrule{8-13}          &       &       &       &       &       &       & \multicolumn{6}{l}{*Optimal solutions obtained by exact algorithms} \\
			&       &       &       &       &       &       & \multicolumn{6}{l}{**Verified as optimal via exact methods} \\
		\end{tabular}%
		\label{tab:op2}%
	\end{sidewaystable}

	\begin{table}
		\centering
		\scriptsize % 1. 使用小号字体 (如果还觉得大，可以改成 \tiny)
		\setlength{\tabcolsep}{4pt} % 2. 极大地减小列间距 (默认是 6pt)
		\renewcommand{\arraystretch}{1}
		\caption{New best known solution on TOP instance.}
		\begin{tabular}{llrrr c llrrr}
			\cmidrule{1-5}\cmidrule{7-11}    instance & B\&C  & \multicolumn{1}{l}{HALNS} & \multicolumn{1}{l}{LNS} & \multicolumn{1}{l}{BKs} &       & instance & \multicolumn{1}{l}{B\&C} & \multicolumn{1}{l}{HALNS} & LNS   & \multicolumn{1}{l}{BKs} \\
			\cmidrule{1-5}\cmidrule{7-11}    p4.2.c & $\checkmark$ &       &       & 452   &       & p7.2.k & \multicolumn{1}{l}{$\checkmark$} &       &       & 705 \\
			p4.2.d & $\checkmark$ &       &       & 531   &       & p7.2.l & \multicolumn{1}{l}{$\checkmark$} &       &       & 767 \\
			p4.2.e & $\checkmark$ &       &       & 618   &       & p7.2.m & \multicolumn{1}{l}{$\checkmark$} &       &       & 827 \\
			p4.2.g & $\checkmark$ &       &       & 757   &       & p7.2.n & \multicolumn{1}{l}{$\checkmark$} &       &       & 888 \\
			p4.2.n & $\checkmark$ &       &       & 1174  &       & p7.2.o & \multicolumn{1}{l}{$\checkmark$} &       &       & 945 \\
			p4.2.o & $\checkmark$ &       &       & 1218  &       & p7.2.p & \multicolumn{1}{l}{$\checkmark$} &       &       & 1002 \\
			p4.3.d & $\checkmark$ &       &       & 335   &       & p7.2.r & \multicolumn{1}{l}{$\checkmark$} &       &       & 1094 \\
			p4.3.e & $\checkmark$ &       &       & 468   &       & p7.2.s & \multicolumn{1}{l}{$\checkmark$} &       &       & 1136 \\
			p4.3.f & $\checkmark$ &       &       & 579   &       & p7.3.o & \multicolumn{1}{l}{$\checkmark$} &       &       & 874 \\
			p4.3.j & $\checkmark$ &       &       & 861   &       & p7.3.p & \multicolumn{1}{l}{$\checkmark$} &       &       & 929 \\
			p4.3.n & $\checkmark$ &       &       & 1121  &       & rd400\_gen2\_m3 &       & \multicolumn{1}{l}{$\checkmark$} &       & 12646 \\
			p4.4.g & $\checkmark$ &       &       & 461   &       & gr229\_gen2\_m4 &       &       & $\checkmark$ & 11359 \\
			p5.3.r & $\checkmark$ &       &       & 1125  &       & gr229\_gen3\_m4 &       &       & $\checkmark$ & 7660 \\
			p7.2.j & $\checkmark$ &       &       & 646   &       & rd400\_gen1\_m2 &       &       & $\checkmark$ & 233 \\
			\cmidrule{1-5}\cmidrule{7-11}
		\end{tabular}
		\label{tab:nbks_top}
	\end{table}
	
	\begin{sidewaystable}
		\centering
		\scriptsize % 1. 使用小号字体 (如果还觉得大，可以改成 \tiny)
		\setlength{\tabcolsep}{1pt} % 2. 极大地减小列间距 (默认是 6pt)
		\renewcommand{\arraystretch}{1}
		\caption{The new best known solution on the TOPTW instance when the number of routes is 1 or 2.}
		\begin{tabular}{clcccccrrllllllll}
			\cmidrule{1-8}\cmidrule{10-17}    \multicolumn{1}{l}{number of route} & instance & \multicolumn{1}{l}{ILS/SAILS} & \multicolumn{1}{l}{CP} & \multicolumn{1}{l}{MS-ILS} & \multicolumn{1}{l}{ADOPT} & \multicolumn{1}{l}{ES} & \multicolumn{1}{l}{BKs} &       & number of route & instance & ILS/SAILS   & CP & MS-ILS & ADOPT & ES & BKs \\
			\cmidrule{1-8}\cmidrule{10-17}    \multirow{17}[2]{*}{1} & r203  & $\checkmark$ &       &       &       &       & 1028  &       & \multicolumn{1}{c}{\multirow{21}[4]{*}{2}} & r209  &       &       & \multicolumn{1}{c}{$\checkmark$} &       &       & \multicolumn{1}{r}{1423} \\
			& r204a & $\checkmark$ &       &       &       &       & 1093  &       &       & r210  &       &       & \multicolumn{1}{c}{$\checkmark$} &       &       & \multicolumn{1}{r}{1438} \\
			& r206a & $\checkmark$ &       &       &       &       & 1032  &       &       & rc201 &       &       & \multicolumn{1}{c}{$\checkmark$} &       &       & \multicolumn{1}{r}{1386} \\
			& r207  &       &       & $\checkmark$ &       &       & 1078  &       &       & rc201a & \multicolumn{1}{c}{$\checkmark$} &       &       &       &       & \multicolumn{1}{r}{1385} \\
			& r208  & $\checkmark$ &       &       &       &       & 1118  &       &       & rc202 &       &       & \multicolumn{1}{c}{$\checkmark$} &       &       & \multicolumn{1}{r}{1523} \\
			& r209  &       &       & $\checkmark$ &       &       & 962   &       &       & rc203 &       &       & \multicolumn{1}{c}{$\checkmark$} &       &       & \multicolumn{1}{r}{1640} \\
			& r210  &       &       & $\checkmark$ &       &       & 1002  &       &       & rc204 &       &       & \multicolumn{1}{c}{$\checkmark$} &       &       & \multicolumn{1}{r}{1718} \\
			& r211  & $\checkmark$ &       &       &       &       & 1051  &       &       & rc206 & \multicolumn{1}{c}{$\checkmark$} &       &       &       &       & \multicolumn{1}{r}{1552} \\
			& rc202a & $\checkmark$ &       &       &       &       & 938   &       &       & rc207 &       &       & \multicolumn{1}{c}{$\checkmark$} &       &       & \multicolumn{1}{r}{1609} \\
			& rc204 &       &       & $\checkmark$ &       &       & 1143  &       &       & rc208 &       &       & \multicolumn{1}{c}{$\checkmark$} &       &       & \multicolumn{1}{r}{1705} \\
			& rc206a & $\checkmark$ &       &       &       &       & 899   &       &       & pr04  &       &       & \multicolumn{1}{c}{$\checkmark$} &       &       & \multicolumn{1}{r}{928} \\
			& rc208 &       &       & $\checkmark$ &       &       & 1058  &       &       & pr04a & \multicolumn{1}{c}{$\checkmark$} &       &       &       &       & \multicolumn{1}{r}{926} \\
			& rc208a & $\checkmark$ &       &       &       &       & 1057  &       &       & pr05  &       &       & \multicolumn{1}{c}{$\checkmark$} &       &       & \multicolumn{1}{r}{1103} \\
			& pr01  &       & $\checkmark$ &       &       &       & \multicolumn{1}{l}{308*} &       &       & pr09a & \multicolumn{1}{c}{$\checkmark$} &       &       &       &       & \multicolumn{1}{r}{909} \\
			& pr13  & $\checkmark$ &       &       &       &       & 467   &       &       & pr10  &       &       & \multicolumn{1}{c}{$\checkmark$} &       &       & \multicolumn{1}{r}{1145} \\
			& pr15  & $\checkmark$ &       &       &       &       & 708   &       &       & pr13  &       &       &       &       & \multicolumn{1}{c}{$\checkmark$} & \multicolumn{1}{r}{845} \\
			& r107a & $\checkmark$ &       &       &       &       & 538   &       &       & pr13a & \multicolumn{1}{c}{$\checkmark$} &       &       &       &       & \multicolumn{1}{r}{843} \\
			\cmidrule{2-8}    \multirow{5}[4]{*}{2} & r201  &       &       & $\checkmark$ &       &       & 1260  &       &       & pr15  &       &       & \multicolumn{1}{c}{$\checkmark$} &       &       & \multicolumn{1}{r}{1238} \\
			& r202  &       &       & $\checkmark$ &       &       & 1353  &       &       & pr18  &       &       & \multicolumn{1}{c}{$\checkmark$} &       &       & \multicolumn{1}{r}{955} \\
			& r203  &       &       & $\checkmark$ &       &       & \multicolumn{1}{l}{1431*} &       &       & pr19  &       &       &       & \multicolumn{1}{c}{$\checkmark$} &       & \multicolumn{1}{r}{1041} \\
			& r205  &       &       & $\checkmark$ &       &       & \multicolumn{1}{l}{1402*} &       &       & pr20  &       &       & \multicolumn{1}{c}{$\checkmark$} &       &       & \multicolumn{1}{r}{1251} \\
			\cmidrule{10-17}          & r206  &       &       & $\checkmark$ &       &       & 1452  &       & \multicolumn{8}{l}{*Verified as optimal via exact methods} \\
			\cmidrule{1-8}    \end{tabular}%
		\label{tab:toptw1}%
	\end{sidewaystable}
	
	\begin{sidewaystable}
		\centering
		\scriptsize % 1. 使用小号字体 (如果还觉得大，可以改成 \tiny)
		\setlength{\tabcolsep}{1pt} % 2. 极大地减小列间距 (默认是 6pt)
		\renewcommand{\arraystretch}{1}
		\caption{The new best known solution on the TOPTW instance when the number of routes is 3 or 4.}
		\begin{tabular}{clccccrrlllllll}
			\cmidrule{1-7}\cmidrule{9-15}    \multicolumn{1}{l}{number of routes} & instance & \multicolumn{1}{l}{ILS/SAILS} & \multicolumn{1}{l}{MS-ILS} & \multicolumn{1}{l}{CP} & \multicolumn{1}{l}{ES} & \multicolumn{1}{l}{BKs} &       & number of routes & instance & ILS/SAILS   & MS-ILS & CP    & ES    & BKs \\
			\cmidrule{1-7}\cmidrule{9-15}    \multirow{16}[2]{*}{3} & r104a & $\checkmark$ &       &       &       & 778   &       & \multicolumn{1}{c}{\multirow{19}[4]{*}{4}} & rc104 & $\checkmark$ &       &       &       & \multicolumn{1}{r}{1065} \\
			& r201  &       & $\checkmark$ &       &       & 1450  &       &       & rc107 & $\checkmark$ &       &       &       & \multicolumn{1}{r}{987} \\
			& rc104a & $\checkmark$ &       &       &       & 835   &       &       & pr02  &       & $\checkmark$ &       &       & \multicolumn{1}{r}{1083} \\
			& pr02  &       & $\checkmark$ &       &       & 947   &       &       & pr03  &       & $\checkmark$ &       &       & \multicolumn{1}{r}{1248} \\
			& pr03  &       & $\checkmark$ &       &       & 1014  &       &       & pr04  &       & $\checkmark$ &       &       & \multicolumn{1}{r}{1595} \\
			& pr04  &       & $\checkmark$ &       &       & 1298  &       &       & pr05  &       & $\checkmark$ &       &       & \multicolumn{1}{r}{1859} \\
			& pr05  &       & $\checkmark$ &       &       & 1500  &       &       & pr06  &       & $\checkmark$ &       &       & \multicolumn{1}{r}{1898} \\
			& pr06  &       & $\checkmark$ &       &       & 1519  &       &       & pr08  &       &       &       & $\checkmark$ & \multicolumn{1}{r}{1392} \\
			& pr08  &       & $\checkmark$ &       &       & 1142  &       &       & pr09  &       & $\checkmark$ &       &       & \multicolumn{1}{r}{1626} \\
			& pr10  &       & $\checkmark$ &       &       & 1582  &       &       & pr10  &       & $\checkmark$ &       &       & \multicolumn{1}{r}{1965} \\
			& pr13  &       & $\checkmark$ &       &       & 1159  &       &       & pr11  &       &       & $\checkmark$ &       & 657* \\
			& pr14  &       & $\checkmark$ &       &       & 1375  &       &       & pr12  &       & $\checkmark$ &       &       & \multicolumn{1}{r}{1135} \\
			& pr15  &       & $\checkmark$ &       &       & 1694  &       &       & pr13  &       & $\checkmark$ &       &       & \multicolumn{1}{r}{1392} \\
			& pr18  &       & $\checkmark$ &       &       & 1289  &       &       & pr14  &       & $\checkmark$ &       &       & \multicolumn{1}{r}{1688} \\
			& pr19  &       & $\checkmark$ &       &       & 1428  &       &       & pr15  &       & $\checkmark$ &       &       & \multicolumn{1}{r}{2085} \\
			& pr20  &       & $\checkmark$ &       &       & 1722  &       &       & pr17  &       & $\checkmark$ &       &       & \multicolumn{1}{r}{936} \\
			\cmidrule{2-7}    \multirow{4}[4]{*}{4} & r104  & $\checkmark$ &       &       &       & 975   &       &       & pr18  &       & $\checkmark$ &       &       & \multicolumn{1}{r}{1558} \\
			& r108  & $\checkmark$ &       &       &       & 995   &       &       & pr19  &       & $\checkmark$ &       &       & \multicolumn{1}{r}{1780} \\
			& r112  & $\checkmark$ &       &       &       & 974   &       &       & pr20  &       & $\checkmark$ &       &       & \multicolumn{1}{r}{2115} \\
			\cmidrule{9-15}          & rc103a & $\checkmark$ &       &       &       & 975   &       & \multicolumn{7}{l}{*Verified as optimal via exact methods} \\
			\cmidrule{1-7}    \end{tabular}%
		\label{tab:toptw2}%
	\end{sidewaystable}

	\section{Assumptions and types of solution approaches of canonical extensions.}\label{append2}
	
	% Table generated by Excel2LaTeX from sheet 'Sheet4'
	
	\begin{sidewaystable}
		\centering
		\scriptsize % 1. 使用小号字体 (如果还觉得大，可以改成 \tiny)
		\setlength{\tabcolsep}{1pt} % 2. 极大地减小列间距 (默认是 6pt)
		\renewcommand{\arraystretch}{1}
		\caption{Assumptions of canonical extended models.}
		\begin{tabular}{llllrrrrcrr}
			\toprule
			Reference  &  Problem  & Team  & \makecell{Time \\windows} & \multicolumn{1}{l}{\makecell{Stochas-\\tic}} & \multicolumn{1}{l}{\makecell{Time-\\Dependent}} & \multicolumn{1}{l}{\makecell{Compul-\\sory}} & \multicolumn{1}{l}{\makecell{Capacit-\\ated}} & \multicolumn{1}{l}{Clustered} & \multicolumn{1}{l}{Dynamic} & \multicolumn{1}{l}{\makecell{Multi-\\object}} \\
			\midrule
			\cite{bian2018} &  OStchOP  &       &       & \multicolumn{1}{c}{$\checkmark$} &       &       &       &       & \multicolumn{1}{c}{$\checkmark$} &  \\
			\cite{dolinskaya2018} &  AOPStchT  &       &       & \multicolumn{1}{c}{$\checkmark$} &       &       &       &       & \multicolumn{1}{c}{$\checkmark$} &  \\
			\cite{Song2020} &  StchTOP-CC  & \multicolumn{1}{c}{$\checkmark$} &       & \multicolumn{1}{c}{$\checkmark$} &       &       &       &       &       &  \\
			\cite{thayer2021} &  StchOP  &       &       & \multicolumn{1}{c}{$\checkmark$} &       &       &       &       &       &  \\
			\cite{panadero2022} &  StchTOP-PDR  & \multicolumn{1}{c}{$\checkmark$} &       & \multicolumn{1}{c}{$\checkmark$} &       &       &       &       &       &  \\
			\cite{carpin2024} &  StchOP  &       &       & \multicolumn{1}{c}{$\checkmark$} &       &       &       &       &       &  \\
			\cite{verbeeck2017} &  TD-OPTW  &       & \multicolumn{1}{c}{$\checkmark$} &       & \multicolumn{1}{c}{$\checkmark$} &       &       &       &       &  \\
			\cite{yu2019a} &  OPSTP  &       &       &       & \multicolumn{1}{c}{$\checkmark$} &       &       &       &       &  \\
			\cite{peng2019} &  OP-TDPT  &       & \multicolumn{1}{c}{$\checkmark$} &       & \multicolumn{1}{c}{$\checkmark$} &       &       &       &       &  \\
			\cite{peng2020} &  OPTW  &       & \multicolumn{1}{c}{$\checkmark$} &       & \multicolumn{1}{c}{$\checkmark$} &       &       &       &       &  \\
			\cite{peng2025} &  TD-TOPVTW  & \multicolumn{1}{c}{$\checkmark$} & \multicolumn{1}{c}{$\checkmark$} &        & \multicolumn{1}{c}{$\checkmark$} &        &        &        &        &  \\
			\cite{yu2019b} &  TOPTW-TDS  & \multicolumn{1}{c}{$\checkmark$} & \multicolumn{1}{c}{$\checkmark$} &       & \multicolumn{1}{c}{$\checkmark$} &       &       &       &       &  \\
			\cite{yu2021} &  TOP-TV  & \multicolumn{1}{c}{$\checkmark$} &       &       & \multicolumn{1}{c}{$\checkmark$} &       &       &       &       &  \\
			\cite{khodadadian2022} &  TD-OPTW-STP  &       & \multicolumn{1}{c}{$\checkmark$} &       & \multicolumn{1}{c}{$\checkmark$} &       &       &       &       &  \\
			\cite{dasdemir2022} &  MO-OP-TDMC  &       &       &       & \multicolumn{1}{c}{$\checkmark$} &       &       &       &       & \multicolumn{1}{c}{$\checkmark$} \\
			\cite{avraham2023} &  \makecell[l]{DD-TD-OP-\\STW}  &       & \multicolumn{1}{c}{$\checkmark$} & \multicolumn{1}{c}{$\checkmark$} & \multicolumn{1}{c}{$\checkmark$} &       &       &       &       &  \\
			\cite{palomo2017} &  OPMVEC  &       &       &       &       & \multicolumn{1}{c}{$\checkmark$} &       &       &       &  \\
			\cite{lu2018} &  OPMVEC  &       &       &       &       & \multicolumn{1}{c}{$\checkmark$} &       &       &       &  \\
			\cite{lin2017} &  TOPTW-MV  & \multicolumn{1}{c}{$\checkmark$} & \multicolumn{1}{c}{$\checkmark$} &       &       & \multicolumn{1}{c}{$\checkmark$} &       &       &       &  \\
			\cite{li2022} &  DTOP-MV  & \multicolumn{1}{c}{$\checkmark$} &       &       &       & \multicolumn{1}{c}{$\checkmark$} &       &       &       &  \\
			\cite{fang2023} &  TOPMV  & \multicolumn{1}{c}{$\checkmark$} &       &       &       & \multicolumn{1}{c}{$\checkmark$} &       &       &       &  \\
			\cite{assuncao2019} &  SteinerTOP  & \multicolumn{1}{c}{$\checkmark$} &       &       &       & \multicolumn{1}{c}{$\checkmark$} &       &       &       &  \\
			\cite{assuncao2021} &  SteinerTOP  & \multicolumn{1}{c}{$\checkmark$} &       &       &       & \multicolumn{1}{c}{$\checkmark$} &       &       &       &  \\
			\cite{perez2025} &  OPMVC  &        &        &        &        & \multicolumn{1}{c}{$\checkmark$} &        &        &        &  \\
			\cite{guastalla2025} &  TOP-ST-MIN  & \multicolumn{1}{c}{$\checkmark$} &        &        &        & \multicolumn{1}{c}{$\checkmark$} &        &        &        &  \\
			\cite{wang2025} &  OPMV-SP  &        &        &        &        & \multicolumn{1}{c}{$\checkmark$} &        &        &        &  \\
			\cite{park2017} &  CTOPTW  & \multicolumn{1}{c}{$\checkmark$} & \multicolumn{1}{c}{$\checkmark$} &       &       &       & \multicolumn{1}{c}{$\checkmark$} &       &       &  \\
			\cite{benSaid2019} &  CTOP  & \multicolumn{1}{c}{$\checkmark$} &       &       &       &       & \multicolumn{1}{c}{$\checkmark$} &       &       &  \\
			\cite{shiri2024} &  CTOP  & \multicolumn{1}{c}{$\checkmark$} &       &       &       &       & \multicolumn{1}{c}{$\checkmark$} &       & \multicolumn{1}{c}{$\checkmark$} &  \\
			\cite{hammami2024} &  CTOP  & \multicolumn{1}{c}{$\checkmark$} &       &       &       &       & \multicolumn{1}{c}{$\checkmark$} &       &       &  \\
			\cite{yahiaoui2019} &  CluTOP  & \multicolumn{1}{c}{$\checkmark$} &       &       &       &       &       & $\checkmark$ &       &  \\
			\cite{wu2024} &  CluOP  &       &       &       &       &       &       & $\checkmark$ &       &  \\
			\cite{he2024} &  CluTOP  & \multicolumn{1}{c}{$\checkmark$} &       &       &       &       &       & $\checkmark$ &       &  \\
			\cite{montemanni2024} &  CluOP  &       &       &       &       &       &       & $\checkmark$ &       &  \\
			\cite{almeida2025} & CluOPS & & & & & & & $\checkmark$ & & \\
			\bottomrule
		\end{tabular}%
		\label{tab:12}%
	\end{sidewaystable}%
	
	\begin{sidewaystable}
		\centering
		\scriptsize % 1. 使用小号字体 (如果还觉得大，可以改成 \tiny)
		\setlength{\tabcolsep}{1pt} % 2. 极大地减小列间距 (默认是 6pt)
		\renewcommand{\arraystretch}{1}
		\caption{Types of solution approaches for canonical extended models.}
		\begin{tabular}{llcccccc}
			\toprule
			Reference  &  Algorithm & \multicolumn{1}{l}{Exact} & Heuristic  & Metaheuristic  & Matheuristic  & Approximate & \makecell{Learning-\\Based}  \\
			\midrule
			\cite{bian2018} &  SMPA &       & \multicolumn{1}{c}{$\checkmark$} & \multicolumn{1}{c}{$\checkmark$} &       &       &  \\
			\cite{dolinskaya2018} &  VNS  &       & \multicolumn{1}{c}{$\checkmark$} & \multicolumn{1}{c}{$\checkmark$} &       &       &  \\
			\cite{Song2020} &  \makecell[l]{Anticipatory Consistent Customer \\Assignments} &       & \multicolumn{1}{c}{$\checkmark$} &       & \multicolumn{1}{c}{$\checkmark$} &       &  \\
			\cite{thayer2021} &  CMDP-based Path Tree &       & \multicolumn{1}{c}{$\checkmark$} &       &       &       &  \\
			\cite{panadero2022} &  Simheuristics &       &       & \multicolumn{1}{c}{$\checkmark$} &       &       &  \\
			\cite{carpin2024} &  MCTS-based Approach &       & \multicolumn{1}{c}{$\checkmark$} &       &       &       &  \\
			\cite{verbeeck2017} &  Fast Ant Colony System &       &       &   \multicolumn{1}{c}{$\checkmark$}    &  &       &  \\
			\cite{yu2019a} &  Two-phase Matheuristic &       &       &       & \multicolumn{1}{c}{$\checkmark$} &       &  \\
			\cite{peng2019} &  BDP-ILS &       &       & \multicolumn{1}{c}{$\checkmark$} &       &       &  \\
			\cite{peng2020} &  ADP-DSSR & $\checkmark$ &       &       &       &       &  \\
			\cite{peng2025} &  Branch-and-Cut-and-Price & $\checkmark$ &        &        &        &        &  \\
			\cite{yu2019b} &  Hybrid ABC &       &       & \multicolumn{1}{c}{$\checkmark$} &       &       &  \\
			\cite{yu2021} & Benders B\&C/ILS-MCS & $\checkmark$ &       &     $\checkmark$  &       &       &  \\
			\cite{khodadadian2022} &  Variable Neighborhood Search &       &       & \multicolumn{1}{c}{$\checkmark$} &       &       &  \\
			\cite{dasdemir2022} &  Hybrid Algorithm &       &       &       & \multicolumn{1}{c}{$\checkmark$} &       &  \\
			\cite{avraham2023} &  Branch-and-Bound & $\checkmark$ &       &       &       &       &  \\
			\cite{palomo2017} &  Hybrid GRASP-VNS &       &       & \multicolumn{1}{c}{$\checkmark$} &       &       &  \\
			\cite{lu2018} &  Memetic Algorithm &       &       & \multicolumn{1}{c}{$\checkmark$} &       &       &  \\
			\cite{lin2017} &  Multi-start Simulated Annealing &       &       & \multicolumn{1}{c}{$\checkmark$} &       &       &  \\
			\cite{li2022} &  IoT-Agile Optimization &       &       & \multicolumn{1}{c}{$\checkmark$} &       &       &  \\
			\cite{fang2023} &  NNH with HRGCN &       &       & \multicolumn{1}{c}{$\checkmark$} &       &       & \multicolumn{1}{c}{$\checkmark$} \\
			\cite{assuncao2019} &  Cutting-plane & $\checkmark$ &       &       &       &       &  \\
			\cite{assuncao2021} &  LNS + FP matheuristic &       &       & \multicolumn{1}{c}{$\checkmark$} &   $\checkmark$    &       &  \\
			\cite{perez2025} &  Branch-and-Cut & $\checkmark$ &        &        &        &        &  \\
			\cite{guastalla2025} &  Cutting-Plane Algorithm & $\checkmark$ &        &        &        &        &  \\
			\cite{wang2025} &  Cutting-Plane Algorithm & $\checkmark$ &        &        &        &        &  \\
			\cite{park2017} &  Branch-and-price & $\checkmark$ &       &       &       &       &  \\
			\cite{benSaid2019} &  Variable Space Search &       &       & \multicolumn{1}{c}{$\checkmark$} &       &       &  \\
			\cite{shiri2024} &  Online algorithms &       &       &       &       & \multicolumn{1}{c}{$\checkmark$} &  \\
			\cite{hammami2024} &  Hybrid ALNS &       &       & \multicolumn{1}{c}{$\checkmark$} & \multicolumn{1}{c}{$\checkmark$} &       &  \\
			\cite{yahiaoui2019} &  Cutting plane method with ALNS & $\checkmark$ &       &    $\checkmark$   &  &       &  \\
			\cite{wu2024} &  Hybrid evolutionary algorithm &       &       & \multicolumn{1}{c}{$\checkmark$} &       &       & \multicolumn{1}{c}{$\checkmark$} \\
			\cite{he2024} &  Multi-level memetic search &       &       & \multicolumn{1}{c}{$\checkmark$} &       &       &  \\
			\cite{montemanni2024} &  Constraint programming model & $\checkmark$ &       &       &       &       &  \\
			\cite{almeida2025} & Branch-and-Cut \& Tabu Search & $\checkmark$ & & \multicolumn{1}{c}{$\checkmark$} & & & \\
			\bottomrule
		\end{tabular}%
		\label{tab:13}%
	\end{sidewaystable}%

	\section{Assumptions and types of solution approaches of emerging new models.}\label{append3}

	\begin{sidewaystable}
		\centering
		\scriptsize % 1. 使用小号字体 (如果还觉得大，可以改成 \tiny)
		\setlength{\tabcolsep}{1pt} % 2. 极大地减小列间距 (默认是 6pt)
		\renewcommand{\arraystretch}{1}
		\caption{Assumptions of emerging new models.}
		\begin{tabular}{llcccccrcc}
			\toprule
			Reference  &  Problem  & Team  & \makecell{Time \\windows} & Set   & \multicolumn{1}{l}{\makecell{Compul\\sory}} & \multicolumn{1}{l}{Dynamic} & \multicolumn{1}{l}{Dubins} & \multicolumn{1}{l}{\makecell{Probabili-\\stic}} & \multicolumn{1}{l}{\makecell{Multi-\\object}} \\
			\midrule
			\cite{archetti2018} &  SOP  &       &       & \multicolumn{1}{c}{$\checkmark$} &       &       &       &       &  \\
			\cite{penicka2019} &  SOP  &       &       & \multicolumn{1}{c}{$\checkmark$} &       &       &   \multicolumn{1}{c}{$\checkmark$}    &       &  \\
			\cite{dutta2020} &  MOSOP  &       &       & \multicolumn{1}{c}{$\checkmark$} &       &       &       &       &  \\
			\cite{carrabs2021} &  SOP  &       &       & \multicolumn{1}{c}{$\checkmark$} &       &       &       &       &  \\
			\cite{dontas2023} &  SOP  &       &       & \multicolumn{1}{c}{$\checkmark$} &       &       &       &       &  \\
			\cite{archetti2024} &  SOP  &       &       & \multicolumn{1}{c}{$\checkmark$} &       &       &       &       &  \\
			\cite{lu2024} &  SOP  &       &       & \multicolumn{1}{c}{$\checkmark$} &       &       &       &       &  \\
			\cite{lin2024} &  SOPMV  &       &       & \multicolumn{1}{c}{$\checkmark$} & \multicolumn{1}{c}{$\checkmark$} &       &       &       &  \\
			\cite{yu2024} &  STOPTW  & \multicolumn{1}{c}{$\checkmark$} & \multicolumn{1}{c}{$\checkmark$} & \multicolumn{1}{c}{$\checkmark$} &       &       &       &       &  \\
			\cite{nguyen2025} & STOP & $\checkmark$ & & $\checkmark$ & & & & & \\
			\cite{penicka2017} &  DOP  &       &       &       &       &       & \multicolumn{1}{c}{$\checkmark$} &       &  \\
			\cite{tsiogkas2018} &  DCOP  &       &       &       &       &       & \multicolumn{1}{c}{$\checkmark$} &       &  \\
			\cite{macharet2021} &  MEDOP  &       &       &       &       &       & \multicolumn{1}{c}{$\checkmark$} &       & $\checkmark$ \\
			\cite{faigl2020} &  CEDOP  &   $\checkmark$    &       &       &       &       & \multicolumn{1}{c}{$\checkmark$} &       &  \\
			\cite{angelelli2017} &  POP  &       &       &       &       &       &       & $\checkmark$ &  \\
			\cite{chou2020} &  POP  &       &       &       &       &       &       & $\checkmark$ &  \\
			\cite{montemanni2025} & POP & & & & & & & $\checkmark$ & \\
			\cite{herrera2022} &  TOP-PD  & \multicolumn{1}{c}{$\checkmark$} &       &       &       &       &       & $\checkmark$ &  \\
			\cite{angelelli2021} &  DPOP  &       &       &       &       & $\checkmark$ &       & $\checkmark$ &  \\
			\bottomrule
		\end{tabular}%
		\label{tab:18}%
	\end{sidewaystable}%
	
	\begin{sidewaystable}
		\centering
		\scriptsize % 1. 使用小号字体 (如果还觉得大，可以改成 \tiny)
		\setlength{\tabcolsep}{1pt} % 2. 极大地减小列间距 (默认是 6pt)
		\renewcommand{\arraystretch}{1}
		\caption{Types of solution approaches for emerging new models.}
		\begin{tabular}{llcccccc}
			\toprule
			Reference  &  Algorithm & Exact  & \multicolumn{1}{l}{Heuristic} & \makecell{Meta-\\heuristic}  & Matheuristic  & \makecell{Approxi-\\mate} & \makecell{Learning-\\Based}  \\
			\midrule
			\cite{archetti2018} &  MASOP &       &       &       & \multicolumn{1}{c}{$\checkmark$} &       &  \\
			\cite{penicka2019} &  VNS  &       &       & \multicolumn{1}{c}{$\checkmark$} &       &       &  \\
			\cite{dutta2020} &  NSGA-II, SPEA2 &       &       & \multicolumn{1}{c}{$\checkmark$} &       &       &  \\
			\cite{carrabs2021} &  BRKGA &       &       & \multicolumn{1}{c}{$\checkmark$} &       &       &  \\
			\cite{dontas2023} &  \makecell[l]{Adaptive Memory \\Matheuristic} &       &       &       & \multicolumn{1}{c}{$\checkmark$} &       &  \\
			\cite{archetti2024} &  Branch-and-Cut & \multicolumn{1}{c}{$\checkmark$} &       &       &       &       &  \\
			\cite{lu2024} &  HEA  &       &       &   \multicolumn{1}{c}{$\checkmark$}    &  &       &  \\
			\cite{lin2024} &  Simulated Annealing &       &       & \multicolumn{1}{c}{$\checkmark$} &       &       &  \\
			\cite{yu2024} &  SARL &       &       & \multicolumn{1}{c}{$\checkmark$} &       &       & \multicolumn{1}{c}{$\checkmark$} \\
			\cite{nguyen2025} & Branch-and-Price \& LNS & $\checkmark$ & & $\checkmark$ & & & \\
			\cite{penicka2017} &  VNS Heuristic &       &       & \multicolumn{1}{c}{$\checkmark$} &       &       &  \\
			\cite{tsiogkas2018} &  Genetic Algorithm &       &       & \multicolumn{1}{c}{$\checkmark$} &       &       &  \\
			\cite{macharet2021} &  Evolutionary Algorithm &       &       & \multicolumn{1}{c}{$\checkmark$} &       &       &  \\
			\cite{faigl2020} &  GSOA Algorithm &       &       &       &       &       & \multicolumn{1}{c}{$\checkmark$} \\
			\cite{angelelli2017} &  B\&C with matheuristics & \multicolumn{1}{c}{$\checkmark$} &       &       & \multicolumn{1}{c}{$\checkmark$} &       &  \\
			\cite{chou2020} &  TS with Monte Carlo &       &       & \multicolumn{1}{c}{$\checkmark$} &       &       &  \\
			\cite{montemanni2025} & Iterative CP-based Model & $\checkmark$ & & & & & \\
			\cite{herrera2022} &  \makecell[l]{Simheuristic with reliability \\analysis} &       &       & \multicolumn{1}{c}{$\checkmark$} &       &       &  \\
			\cite{angelelli2021} &  SAA+greedy &       & $\checkmark$ &       &       &       & $\checkmark$ \\
			\bottomrule
		\end{tabular}%
		\label{tab:19}%
	\end{sidewaystable}%
	
	\begin{sidewaystable}
		\centering
		\scriptsize % 1. 使用小号字体 (如果还觉得大，可以改成 \tiny)
		\setlength{\tabcolsep}{1pt} % 2. 极大地减小列间距 (默认是 6pt)
		\renewcommand{\arraystretch}{1}
		\caption{Assumptions of other emerging new models.}	
		\begin{tabular}{llccccccccc}
			\toprule
			Reference  &  Problem  & \multicolumn{1}{l}{Team} & \multicolumn{1}{l}{\makecell{Time \\windows}} & \multicolumn{1}{l}{Dynamic} & \multicolumn{1}{l}{Uncertain} & \multicolumn{1}{l}{with Drone} & \multicolumn{1}{l}{Drone routing} & \multicolumn{1}{l}{\makecell{Multi-\\Period}} & \multicolumn{1}{l}{\makecell{Multi-\\Profit}} & \multicolumn{1}{l}{Robust} \\
			\midrule
			\cite{wang2023} & DynOP &       &       & $\checkmark$ &       &       &       &       &       &  \\
			\cite{saeedvand2020} & TOPTWR & $\checkmark$ & \multicolumn{1}{c}{$\checkmark$} & $\checkmark$ &       &       &       &       &       &  \\
			\cite{uguina2024} & DynTOP & $\checkmark$ &       & $\checkmark$ &       &       &       &       &       &  \\
			\cite{li2024} & DyOP-rd &       &       & $\checkmark$ &       &       &       &       &       &  \\
			\cite{granda2025} & TOPVTW & $\checkmark$ & \multicolumn{1}{c}{$\checkmark$} & $\checkmark$ &        &        &        &        &        &  \\
			\cite{kirac2025} & DTOP & $\checkmark$ &        & $\checkmark$ &        &        &        &        &        &  \\
			\cite{paradiso2025} & DTOP-TW & $\checkmark$ & \multicolumn{1}{c}{$\checkmark$} & $\checkmark$ & \multicolumn{1}{c}{$\checkmark$} &        &        &        &        &  \\
			\cite{wang2019a} & UOP   &       &       &       & \multicolumn{1}{c}{$\checkmark$} &       &       &       &       &  \\
			\cite{qian2025} & UDyOP &        &        & $\checkmark$ & \multicolumn{1}{c}{$\checkmark$} &        &        &        &        &  \\
			\cite{jin2019} & TOPSRST & $\checkmark$ &       &       & \multicolumn{1}{c}{$\checkmark$} &       &       &       &       &  \\
			\cite{wang2019b} & UTOPTW & $\checkmark$ & \multicolumn{1}{c}{$\checkmark$} &       & \multicolumn{1}{c}{$\checkmark$} &       &       &       &       &  \\
			\cite{yu2022} & RTOP-DP & $\checkmark$ & \multicolumn{1}{c}{$\checkmark$} &       &       &       &       &       &       & $\checkmark$ \\
			\cite{zhang2023} & RAOP  & $\checkmark$ &       &       &       &       &       &       &       & $\checkmark$ \\
			\cite{shi2023} & RMOP  & $\checkmark$ &       & $\checkmark$ &       &       &       &       &       & $\checkmark$ \\
			
			\cite{bayliss2020} & TOP-D & $\checkmark$ &       &       &       &       & \multicolumn{1}{c}{$\checkmark$} &       &       &  \\
			\cite{sundar2022} & TOP-D & $\checkmark$ &       &       &       &       & \multicolumn{1}{c}{$\checkmark$} &       &       &  \\
			\cite{asghar2024} & MML-TOP & $\checkmark$ &        &        &        & \multicolumn{1}{c}{$\checkmark$} &        &        &        &  \\
			\cite{qian2024} & CEOP &        &        &        &        &  &   $\checkmark$     &        &        &  \\
			\cite{wan2024} & TOP-TV & $\checkmark$ &        &        &        & \multicolumn{1}{c}{$\checkmark$} &        &        & \multicolumn{1}{c}{$\checkmark$} &  \\
			\cite{zhang2020} & MP-OP &       & \multicolumn{1}{c}{$\checkmark$} & $\checkmark$ &       &       &       & \multicolumn{1}{c}{$\checkmark$} &       &  \\
			\cite{wang2024} & MPOP &        &        &        &        &        &        & \multicolumn{1}{c}{$\checkmark$} &        &  \\
			\cite{vidigal2023} & MP-OP &       & \multicolumn{1}{c}{$\checkmark$} &       &       &       &       & \multicolumn{1}{c}{$\checkmark$} &       &  \\
			\cite{kim2020} & MP-OP &       &       &       &       &       &       &       & \multicolumn{1}{c}{$\checkmark$} &  \\
			\cite{kim2022} & MP-OP &       &       &       &       &       &       &       & \multicolumn{1}{c}{$\checkmark$} &  \\
			\cite{morandi2024} & OP-mD &  &       &       &       & \multicolumn{1}{c}{$\checkmark$} &       &       &       &  \\
			
			\bottomrule
		\end{tabular}%
		\label{tab:20}%
	\end{sidewaystable}%
	
	\begin{sidewaystable}
		\centering
		\scriptsize % 1. 使用小号字体 (如果还觉得大，可以改成 \tiny)
		\setlength{\tabcolsep}{1pt} % 2. 极大地减小列间距 (默认是 6pt)
		\renewcommand{\arraystretch}{1}
		\caption{Types of solution approaches for emerging new models.}
		\begin{tabular}{llcccccc}
			\toprule
			Reference  &  Algorithm & \multicolumn{1}{l}{Exact } & \multicolumn{1}{l}{Heuristic } & \multicolumn{1}{l}{Metaheuristic } & \multicolumn{1}{l}{Matheuristic } & \multicolumn{1}{l}{Approximate} & \multicolumn{1}{l}{Learning-Based } \\
			\midrule
			\cite{wang2023} & IRH, HSATSRO, TRO &       &       & \multicolumn{1}{c}{$\checkmark$} &       &       &  \\
			\cite{saeedvand2020} & RL    &       &       & \multicolumn{1}{c}{$\checkmark$} &       &       & \multicolumn{1}{c}{$\checkmark$} \\
			\cite{uguina2024} & Learnheuristic with TS &       & \multicolumn{1}{c}{$\checkmark$} &       &       &       & \multicolumn{1}{c}{$\checkmark$} \\
			\cite{li2024} & VFA-CF and VFA-2S &       &       &       &       &       & \multicolumn{1}{c}{$\checkmark$} \\
			\cite{granda2025} & Greedy Constructive + MIP &        &        &        & \multicolumn{1}{c}{$\checkmark$} &        &  \\
			\cite{kirac2025} & MPA + ALNS &        &        & \multicolumn{1}{c}{$\checkmark$} &        &        &  \\
			\cite{paradiso2025} & Stochastic Lookahead + CG &        &        &        & \multicolumn{1}{c}{$\checkmark$} &        &  \\
			\cite{wang2019a} & one-dimensional ratio search algorithm &       &       &       & \multicolumn{1}{c}{$\checkmark$} &       &  \\
			\cite{qian2025} & ADAPT (Bayesian + IACS) &        &        & \multicolumn{1}{c}{$\checkmark$} &        &        & \multicolumn{1}{c}{$\checkmark$} \\
			\cite{jin2019} & Monte Carlo + heuristics &       &       & \multicolumn{1}{c}{$\checkmark$} &       &       &  \\
			\cite{wang2019b} & Iterated local search &       &  &   $\checkmark$    &       &       &  \\
			\cite{yu2022} & B\&P/TS    & \multicolumn{1}{c}{$\checkmark$} &       & \multicolumn{1}{c}{$\checkmark$} &       &       &  \\
			\cite{zhang2023} & Graph transform + B\&P & \multicolumn{1}{c}{$\checkmark$} &       &       &       &       &  \\
			\cite{shi2023} & Online: MCTS &       &   $\checkmark$    &       &       & $\checkmark$ &  \\
			
			\cite{bayliss2020} & Learnheuristic &       &       & \multicolumn{1}{c}{$\checkmark$} &       &       & \multicolumn{1}{c}{$\checkmark$} \\
			\cite{sundar2022} & DSSR + branch-and-price & \multicolumn{1}{c}{$\checkmark$} &       &       &       &       &  \\
			\cite{asghar2024} & Orienteering-based Heuristic &        & \multicolumn{1}{c}{$\checkmark$} &        &        & \multicolumn{1}{c}{$\checkmark$} &  \\
			\cite{qian2024} & RSZD + Hybrid ACS &        &        & \multicolumn{1}{c}{$\checkmark$} &        &        &  \\
			\cite{wan2024} & Attention-based DRL &        &        &        &        &        & \multicolumn{1}{c}{$\checkmark$} \\		
			\cite{zhang2020} & MDP model &       & \multicolumn{1}{c}{$\checkmark$} &       & \multicolumn{1}{c}{$\checkmark$} &    \checkmark   &  \\
			\cite{wang2024} & Lagrangian Relaxation &        &        &        & \multicolumn{1}{c}{$\checkmark$} &        &  \\
			\cite{vidigal2023} & ILS   &       &       & \multicolumn{1}{c}{$\checkmark$} &       &       &  \\
			\cite{kim2020} & SA + VTOA &       &       & \multicolumn{1}{c}{$\checkmark$} & \multicolumn{1}{c}{$\checkmark$} &       &  \\
			\cite{kim2022} & Hybrid DP & \multicolumn{1}{c}{$\checkmark$} &       &       &       &       &  \\
			\cite{morandi2024} & B\&C  & \multicolumn{1}{c}{$\checkmark$} &       &       &       &       &  \\
			\bottomrule
		\end{tabular}%
		\label{tab:21}%
	\end{sidewaystable}%

\end{document}